\documentclass[journal]{IEEEtran}
\usepackage{cite}
\ifCLASSINFOpdf
   \usepackage[pdftex]{graphicx}
   \DeclareGraphicsExtensions{.pdf,.jpeg,.png}
\else
   \usepackage[dvips]{graphicx}
   \DeclareGraphicsExtensions{.eps}
\fi
\usepackage{color}

\usepackage{amsmath,amssymb}
\usepackage{sobolev}
\usepackage{stmaryrd}
\usepackage{array}
\usepackage{url}

\newcommand{\normdev}{\partial_{\mathbf{n}}}
\renewcommand{\j}{\mathbf{j}}
\newcommand{\n}{\mathbf{n}}
\newcommand{\x}{\mathbf{x}}

\newcommand{\R}{\mathbb{R}}
\renewcommand{\div}{\operatorname{div}}
\newcommand{\Hdivstar}{H_{0}(\div,\Omega)}
\newcommand{\dx}{\,\mathrm{d}\mathit{x}}

\newcommand{\supp}{\operatorname{supp}}
\newcommand{\sign}{\operatorname{sign}}
\renewcommand{\forall}{\text{ for all }}
\renewcommand{\ker}{\operatorname{ker}}

\newtheorem{defi}{Definition}
\newtheorem{rema}[defi]{Remark}

\makeatletter
\newcommand{\pushright}[1]{\ifmeasuring@#1\else\omit\hfill$\displaystyle#1$\fi\ignorespaces}
\newcommand{\pushleft}[1]{\ifmeasuring@#1\else\omit$\displaystyle#1$\hfill\fi\ignorespaces}
\makeatother

\def\uscore{\protect\rule{.2em}{.2pt}}
\def\Uscore{\protect\rule{.35em}{.2pt}}

\newcommand{\segres}[2]{\textit{seg\uscore{#1}\Uscore{}res\uscore{#2}}}
\newcommand{\segresR}[3]{\textit{seg\uscore{#1}\Uscore{}res\uscore{#2}\Uscore{}r{#3}}}
\newcommand{\ci}[3]{\textit{{#1}C\Uscore{#2}\Uscore{#3}}}

\newcommand{\JP}[1]{[#1]} % point-wise jump
\newcommand{\av}[1]{\{ #1 \}}

\usepackage[final]{changes}

\usepackage{tikz}
\usepackage{textcomp}
\usepackage{hyperref}
\newcommand\copyrighttext{%
  \footnotesize \textcopyright 2016 IEEE. Personal use of this material is permitted.
  Permission from IEEE must be obtained for all other uses, in any current or future 
  media, including reprinting/republishing this material for advertising or promotional 
  purposes, creating new collective works, for resale or redistribution to servers or 
  lists, or reuse of any copyrighted component of this work in other works. 
  DOI: \href{https://doi.org/10.1109/TMI.2016.2624634}{10.1109/TMI.2016.2624634}}
\newcommand\copyrightnotice{%
\begin{tikzpicture}[remember picture,overlay]
\node[anchor=south,yshift=0pt] at (current page.south) {\fbox{\parbox{\dimexpr\textwidth-\fboxsep-\fboxrule\relax}{\copyrighttext}}};
\end{tikzpicture}%
}

\begin{document}
\title{A Mixed Finite Element Method \\ to Solve the EEG Forward Problem}

\author{J.~Vorwerk*, C.~Engwer, S.~Pursiainen, and~C.H.~Wolters% 
\thanks{Copyright (c) 2016 IEEE. Personal use of this material is permitted. However, permission to use this material for any other purposes must be obtained from the IEEE by sending a request to pubs-permissions@ieee.org.}
\thanks{\textit{Asterisk indicates corresponding author.}}
\thanks{J.V. is with the Institute for Biomagnetism and Biosignalanalysis, University of M\"unster, Germany, and the Scientific Computing and Imaging (SCI) Institute, University of Utah, Salt Lake City, USA. e-mail: jvorwerk@sci.utah.edu.}% <-this % stops a space
\thanks{C.E. is with the Institute for Applied Mathematics, University of M\"unster, Germany, and the Cluster of Excellence EXC 1003, Cells in Motion, CiM, M\"unster, Germany.}
\thanks{S.P. is with the Department of Mathematics, Tampere University of Technology, Finland, and the Department of Mathematics and System Analysis, Aalto University, Helsinki, Finland.}
\thanks{C.H.W. is with the Institute for Biomagnetism and Biosignalanalysis, University of M\"unster, Germany.}
\thanks{J.V. and C.H.W. were supported by the Priority Program 1665 of the Deutsche Forschungsgemeinschaft
(DFG) (WO1425/5-2, WO1425/7-1) and the EU project ChildBrain (Marie Curie Innovative
Training Networks, grant agreement no. 641652). C.E. was supported by the Cluster of Excellence 1003 of the Deutsche Forschungsgemeinschaft
(DFG EXC 1003 Cells in Motion). S.P. was supported by the Academy of Finland (Centre of Excellence in Inverse Problems Research and Key Project number 305055).}% <-this % stops a space
}


\markboth{IEEE Transactions on Medical Imaging}%
{Vorwerk \MakeLowercase{\textit{et al.}}: A Mixed Finite Element Method to Solve the EEG Forward Problem}



\maketitle

\copyrightnotice

\begin{abstract}
Finite element methods have been shown to achieve high accuracies in numerically solving the EEG forward problem and they enable the realistic modeling of complex geometries and important conductive features such as anisotropic conductivities. To date, most of the presented approaches rely on the same underlying formulation, the continuous Galerkin (CG)-FEM. In this article, a novel approach to solve the EEG forward problem based on a mixed finite element method (Mixed-FEM) is introduced. To obtain the Mixed-FEM formulation, the electric current is introduced as an additional unknown besides the electric potential. As a consequence of this derivation, the Mixed-FEM is, by construction, current preserving, in contrast to the CG-FEM. Consequently, a higher simulation accuracy can be achieved in certain scenarios, e.g., when the diameter of thin insulating structures, such as the skull, is in the range of the mesh resolution.

A theoretical derivation of the Mixed-FEM approach for EEG forward simulations is presented, and the algorithms implemented for solving the resulting equation systems are described. Subsequently, first evaluations in both sphere and realistic head models are presented, and the results are compared to previously introduced CG-FEM approaches. Additional visualizations are shown to illustrate the current preserving property of the Mixed-FEM.

Based on these results, it is concluded that the newly presented Mixed-FEM can at least complement and in some scenarios even outperform the established CG-FEM approaches, which motivates a further evaluation of the Mixed-FEM for applications in bioelectromagnetism.
\end{abstract}


% \fi
\IEEEpeerreviewmaketitle

\pagebreak

\section{Introduction}
\IEEEPARstart{T}{he} EEG forward problem is to simulate the electric potential on the head surface that is generated by a minimal patch of active brain tissue. Its accurate solution is fundamental for precise EEG source analysis. An accurate solution can be achieved via numerical methods that allow to take the realistic head geometry into account. In this context, finite element methods (FEM) achieve high numerical accuracies and enable to realistically model tissue boundaries with complicated shapes, such as the gray matter/CSF interface, and to incorporate tissue conductivity anisotropy. The importance of incorporating these model features for the computation of accurate forward solutions and, in consequence, also for precise source analysis has been shown in multiple studies \cite{aca2013effects,cho2015influence,CHW:Vor2014}.

Different FEM approaches to solve the EEG forward problem have been proposed, e.g., St. Venant, partial integration, Whitney, or subtraction approaches \cite{CHW:Buc97,CHW:Yan91,CHW:Pur2011,CHW:Wol2007e,vor2016,CHW:Wol2007b}. These approaches differ in the way the dipole source is modeled, but the underlying discretization of the continuous partial differential equation (PDE) is the same: a conforming Galerkin-FEM (CG-FEM) with most often linear Ansatz-functions. The necessary discretization of the head volume can be achieved using either tetrahedral or hexahedral head models. Hexahedral models have the advantage that they can be directly generated from voxel-based magnetic resonance images
(MRI), whereas the generation of 
%triangular surfaces of the compartment boundaries as a basis for
 surface-based tetrahedral meshes can be complicated. 
 Therefore, hexahedral meshes are more and more frequently used in praxis \cite{ayd2014combining,CHW:Rul2009} and have, furthermore, recently been positively validated in an animal study \cite{lau2016skull}.

%A drawback when appyling CG-FEM with hexahedral meshes is the possibility of leakage effects, decreasing the accuracy of the forward solution \cite{son2013leakage}. These occur when the thickness of a thin, low-conducting tissue structure, e.g., the skull, separating two compartments with a higher conductivity, e.g., CSF and skin, is in the range of the mesh resolution as it might be the case in temporal areas, where the skull is only a few millimeters thin. In this case, elements of CSF and skin might touch in single vertices and a current flow through these vertices occurs (Figure \ref{fig:current}). Besides manually increasing the thickness of the skull segmentation in the affected areas, this unwanted -- since unphysical -- effect can be avoided by the use of FEM approaches that explicitly control the flow of volume currents, such as discontinuous Galerkin (DG) or mixed methods. These approaches by construction allow current flow only through faces, but not through single vertices. Thereby, a physically reasonable current flow is enforced and leakage effects can be avoided.

In this article, a mixed finite element method (Mixed-FEM) to solve the EEG forward problem is introduced. 
%Besides 
Compared to the CG-FEM, it has the advantage that the current source can be represented in a direct way, whereas either an approximation using electrical monopoles has to be derived or the subtraction approach has to be applied when using the CG-FEM. Furthermore, the Mixed-FEM is current preserving and thereby prevents the effects of the (local) current leakages through the skull that might occur for the CG-FEM \cite{son2013leakage,eng2015subtraction}. Mixed- and CG-FEM are compared in such a leakage scenario in Section \ref{sub:results-leaky}. An accurate simulation of the currents penetrating the skull is important, as the influence of an accurate representation of the skull for accurate forward simulations has been shown \cite{mon2014influence,CHW:Lan2012,CHW:Dan2011}. 
The accuracy of the Mixed-FEM in comparison to CG-FEM approaches and a recently presented approach based on a discontinuous Galerkin (DG) FEM formulation\mbox{\cite{eng2015subtraction}} is evaluated in sphere and realistic head models. It is shown that the Mixed-FEM achieves higher accuracies in solving the EEG forward problem than the CG-FEM for highly eccentric sources in sphere models and than both CG- and DG-FEM in realistic head models.
%short: with leakage effects.

\section{Theory}
\IEEEPARstart{A}{pplying} the quasistatic approximation of Maxwell's equations \cite{ham1993,CHW:Bre2012}, the forward problem of EEG is commonly formulated as a second-order PDE with homogeneous Neumann boundary condition
\begin{subequations}
\label{eq:forward}
\begin{align}
 \nabla \cdot ( \sigma \nabla u ) &= \nabla \cdot \j^p &  &\text{ in } \Omega, \label{eq:forward1}\\
\sigma \normdev u &= 0 & &\text{ on } \partial \Omega = \Gamma. \label{eq:forward2}
\end{align}
\end{subequations}
Here, $u$ denotes the electric potential, $\j^p$ the source current, and $\sigma$ the conductivity distribution in $\Omega$. In \eqref{eq:forward}, the electric current $\j$ is already eliminated as an unknown. For our purpose, we start at the previous step\added{ in the derivation of the quasi\deleted{-}static approximation} and keep the electric current as an unknown. Thus, our starting point is the system of first-order PDEs
\begin{subequations}
\label{eq:forward-system}
\begin{align}
   \j + \sigma \nabla u &= \j^p &\label{eq:forward-system1}\\\
 \nabla \cdot \j &= 0
\quad \text{ in } \Omega, \label{eq:forward-system2}\\
 \langle \j , \n \rangle &= \langle \j^p , \n \rangle \quad \text{ on } \partial \Omega = \Gamma. \label{eq:forward-system3}
\end{align}
\end{subequations}
Since the source current $\j^p$ in general fulfills $\langle \j^p, \n \rangle = 0$ on $\Gamma$, as $\supp \j^p \subset \Omega^\circ$ for physiological reasons (there are no sources in the skin), \eqref{eq:forward-system3} can be simplified to $\langle \j, \n \rangle = 0$ on $\Gamma$.
The Mixed-FEM formulation for the EEG forward problem is now derived from \eqref{eq:forward-system}, instead of discretizing \eqref{eq:forward} as would be done for the CG-FEM. 

\subsection{A (Mixed) Weak Formulation of the EEG Forward Problem}
\label{sub:mixed-weak-formulation}
Due to the vector-valued equation \eqref{eq:forward-system1}, it is necessary to introduce a space of vector-valued test functions to \replaced{be able to derive}{obtain} a weak formulation of \eqref{eq:forward-system}. A natural function space for the current in the mixed formulation is $\Hdiv$:
 \begin{equation}
  \Hdiv = \left\{ \mathbf{q} \in \L2^3 : \nabla \cdot \mathbf{q} \in \L2 \right\}.
 \end{equation}
 Akin to the scalar-valued Sobolev spaces $\H{k}$, this space becomes a Hilbert space with the norm
 \begin{equation}
 \| \mathbf{q} \|_{\Hdiv} = \left( \| \mathbf{q} \|_{\L2^3}^2 + \| \nabla \cdot \mathbf{q} \|_{\L2}^2 \right)^\frac{1}{2}.
\end{equation}
We introduce a subspace $\Hdivstar$ of $\Hdiv$, in which the boundary condition $\langle \j , \n \rangle = 0 \text{ on } \partial \Omega = \Gamma$ is fulfilled by construction:
 \begin{equation}
 \label{def:hdivstar}
  \Hdivstar = \left\{ \mathbf{q} \in \Hdiv : \langle \mathbf{q}, \mathbf{n} \rangle = 0 \text{ on }\partial \Omega \right
  \}.
 \end{equation}
For the scalar-valued equation \eqref{eq:forward-system2}, one can simply choose the space of square-integrable functions $\L2$ as the test space.

Now, we can introduce a weak formulation of \eqref{eq:forward-system}
\begin{subequations}
\label{eq:mixed-weak}
\begin{multline}
\noindent \int_{\Omega} \langle \sigma^{-1} \j , \mathbf{q} \rangle \dx  - \int_\Omega \nabla \cdot \mathbf{q} u \dx = \\ \int_\Omega \langle \sigma^{-1} \j^p , \mathbf{q} \rangle \dx \quad \forall \mathbf{q} \in \Hdivstar,\label{eq:mixed-weak1}
 \end{multline} 
 \begin{equation}
 \int_\Omega \nabla \cdot \j v \dx = 0 \quad \forall v \in \L2. \label{eq:mixed-weak2}
\end{equation}
\end{subequations}
This is the so-called \emph{dual mixed formulation}\mbox{\cite{rob1991,arn1990mixed,bre2012mixed}}. The Neumann boundary condition \eqref{eq:forward-system3}/\eqref{def:hdivstar} is an essential boundary condition in the (dual) mixed formulation and has to be imposed explicitly in solving the discrete problem.
 We define the bilinear forms 
\begin{subequations}
\label{def:mixed}
\begin{align}
 a(\mathbf{p},\mathbf{q}) &= (\sigma^{-1}\mathbf{p},\mathbf{q})_{\L2^3}, \label{def:mixed-a}\\
 b(\mathbf{p},v) &= (\nabla \cdot \mathbf{p},v)_{\L2} \label{def:mixed-b}
\end{align}
and the functional
\begin{equation}
 l(\mathbf{q}) = (\sigma^{-1}\j^p,\mathbf{q})_{\L2^3}
\end{equation}
\end{subequations}
for $\mathbf{p},\mathbf{q} \in \Hdivstar$, $v \in \L2$, $\j^p \in \L2^3$, $\sigma \in \L{\infty}$, $\sigma > 0$. Therefore, to solve \eqref{eq:mixed-weak} is to
\begin{subequations}
 \label{eq:mixed-op}
\begin{align}
  \intertext{\quad find $(u,\mathbf{j}) \in \L2 \times \Hdivstar$, such that} 
  &a(\mathbf{j},\mathbf{q}) + b(\mathbf{q},u) & \kern-1.05em &= l(\mathbf{q}) & &\forall \mathbf{q} \in \Hdivstar, \label{eq:mixed-op1}\\
  &b(\mathbf{j},v) & \kern-1.05em &= 0 & &\forall v \in \L2. \label{eq:mixed-op2}
 \end{align}
\end{subequations}
In this notation, the saddle point structure of problem \eqref{eq:mixed-op} and thus also \eqref{eq:mixed-weak} is recognizable. As a consequence, the existence and uniqueness of a solution cannot be shown using the Lemma of Lax-Milgram. 

Instead, it can be shown that a solution to \eqref{eq:mixed-op} exists, if the operator $a$ is $\Hdivstar$-elliptic on the kernel of $b$, and $b$ fulfills an \emph{inf-sup condition}, which in this case is also called the \emph{LBB condition}, named after the mathematicians \emph{Ladyzhenskaya}, \emph{Babuska}, and \emph{Brezzi}. At this point, we shall note only that these conditions are fulfilled by $a$ and $b$ defined as in \eqref{def:mixed} and thereby the existence of a solution $(u,\j)\in \L2 \times \Hdivstar$ is given. While the solution for $\j$ is unique, $u$ is defined up to an element of $\ker b(\j,v)$, $v \in \L2$\mbox{\cite{rob1991,ber1994mixed}}. Uniqueness for $u$ can be obtained by introducing an additional condition, such as fixing the value of $\int_\Omega u$ or $\int_{\partial \Omega} u$. For a detailed proof and discussion we refer the reader to, e.g., \mbox{\cite{rob1991,bre2012mixed,CHW:Bra2007}}.

\subsection{Mixed Finite Element Method}
Obtaining a numerical solution for \eqref{eq:mixed-weak}/\eqref{eq:mixed-op} necessitates choosing suitable discrete approximations for the test function spaces $\Hdivstar$ and $\L2$. Utilizing a Galerkin approach, these are also the spaces in which the discrete solution $(u_h, \j_h)$ lies.

In order to construct the discrete subspaces, the volume $\Omega$ is subdivided and approximated by a set of simple geometrical objects. In three dimensions, these objects are usually tetrahedra or hexahedra. For the sake of simplicity, any subdivision of $\Omega$ into either tetrahedra or hexahedra is henceforth referred to as triangulation $\mathcal{T} = \{T_1, T_2, T_3, ..., T_m \}$. In this paper, we follow the definition of a triangulation according to\mbox{\cite{CHW:Bra2007}}. We further assume that the triangulation $\mathcal{T}$ is admissible \cite{CHW:Bra2007} and write $\mathcal{T}_h$ instead of $\mathcal{T}$, if each element $T \in \mathcal{T}$ has a diameter of maximally $2h$.

%\begin{defi}
%\label{def:triangulation}
%A decomposition $\mathcal{T}$ of $\Omega$ is called admissible, iff
%\begin{enumerate}
%\item $\bar{\Omega} = \bigcup_i T_i$,
%\item $T_i \neq \emptyset \forall i$,
%\item for $i \neq j$, $T_i \cap T_j = \partial T_i \cap \partial T_j$ and is (depending on the dimension) either a vertex, an edge, or a face, i.e., $\codim (T_i \cap T_j, \Omega ) > 0$.
%\end{enumerate}
%\end{defi}
%\begin{rema}
%\begin{itemize}
%\item We write $\mathcal{T}_h$ instead of $\mathcal{T}$, if each element $T \in \mathcal{T}$ has a diameter of maximally $2h$.
%\item A family of triangulations $\{ \mathcal{T}_h \}$ is called \emph{shape regular}, if there exists $\kappa > 0$, such that for every $T \in \mathcal{T}_h$ with radius $\rho_T$ of the inscribed circle \begin{equation*}
%\rho_T \geq h_T / \kappa
%\end{equation*}
% holds true, where $h_T$ is half the diameter of $T$.
% \item Though any considered triangulation $\mathcal{T}_h$ is assumed to be a decomposition of the head domain $\Omega$ in the following, we sometimes explicitly refer to it as $\mathcal{T}_h(\Omega)$.
%\end{itemize}
%\end{rema}

We can now choose the space $P_0$ of piecewise constant functions on each element as a discrete subspace of $\L2$:
\begin{align*}
P_0(\mathcal{T}_h) &= \left\{ v \in \L2 : v|_T \equiv c_T,\,c_T \in \R \forall T \in \mathcal{T}_h \right\}.
\end{align*}
A basis of this space is given by the set of characteristic functions $\mathbf{1}_T \in \L2$ for each element $T \in \mathcal{T}_h$. We denote this set of $P_0$ basis functions by $S^{P_0}_h = \{\mathbf{1}_T,\,T \in \mathcal{T}_h\}$.

For $\Hdivstar$, we start by defining the space $RT_0$ of the lowest-order Raviart-Thomas elements on a single, regular hexahedron $T$ \cite{ned1980mixed,rav1977mixed}:
\begin{align*}
RT_0(T) =&  \mathcal{P}_{1,0,0}(T) \times \mathcal{P}_{0,1,0}(T) \times \mathcal{P}_{0,0,1}(T),
\intertext{where $\mathcal{P}_{i,j,k}(T)$ denotes the set of polynomial functions defined on $T$ of degrees $i$, $j$, and $k$ in $x_1$, $x_2$, and $x_3$. We expand this definition to a discrete subspace of $\Hdiv$:}
 RT_0(\mathcal{T}_h) =&  \left\{ \mathbf{q} \in \L2^3 : \mathbf{q}|_T \in RT_0(T) \text{ and } \JP{\mathbf{q}}_{\partial T} = 0 \right. \\ & \pushright{  \forall T \in \mathcal{T}_h \big\}}
\\
 =& \left\{ \mathbf{q} \in \Hdiv : \mathbf{q}|_T \in RT_0(T) \forall T \in \mathcal{T}_h \right\}.
\end{align*}

$\JP{\cdot}_\gamma$ indicates the jump of the normal component at a boundary $\gamma$%, i.e., $\JP{\mathbf{q}}_\gamma = \langle \mathbf{q}^+_\gamma, \n_\gamma \rangle - \langle \mathbf{q}^-_\gamma, \n_\gamma \rangle$ with $\mathbf{q}^+_\gamma$ and $\mathbf{q}^-_\gamma$ indicating the interior and exterior limit value of $\mathbf{q}$ at the boundary $\gamma$
. 

%$RT_0$ can be restricted to a subspace $RT_0^*$ of $\Hdivstar$ through the additional condition $\langle \mathbf{q}, \mathbf{n} \rangle = 0$ on $\partial \Omega$.
Using \emph{Fortin's criterion}\mbox{\cite{bre2012mixed,CHW:Bra2007}}, it can be shown that the existence and uniqueness of a solution to \eqref{eq:mixed-op} -- as noted in Section \ref{sub:mixed-weak-formulation} -- are conserved when replacing $\L2$ and $\Hdivstar$ by their discrete approximations $P_0$ and $RT_0$. For details, we refer the reader to\mbox{\cite{bre2012mixed}}.

A basis of the space $RT_0$ can be defined for both tetrahedral and hexahedral elements. We explicitly note only the hexahedral case, which is also used in the numerical evaluations.
%For tetrahedral elements, a finite dimensional basis of $RT_0$ is 
%%(up to a scaling factor)
% defined by the vector-valued functions supported on two adjacent elements $T_1$ and $T_2$ sharing a common face $f_{1,2} =  \bar{{T}}_1 \cap \bar{{T}}_2 \in \mathcal{T}_h$, so that the normal derivative is continuous across this face and zero on all other faces. With the vertices $\mathbf{v}_i$ defined as in Figure \ref{fig:rt0-support} (left), we have
%
%\begin{equation}
%\label{eq:rt0-basis}
%\mathbf{w}_k ( \x ) = \left\{ \begin{aligned}
%\frac{|f_{1,2}|}{3 |T_1|} \frac{\mathbf{v}_4 - \x}{\| \mathbf{v}_4 - \mathbf{v}_1 \|}_2 &\text{ if } \x \in T_1, \\
%\frac{|f_{1,2}|}{3 |T_2|} \frac{\x - \mathbf{v}_1}{\| \mathbf{v}_4 - \mathbf{v}_1 \|}_2 &\text{ if } \x \in T_2, \\
%0\, &\text{ otherwise.}
%\end{aligned} \right.
%\end{equation}
%
%$|T_1|$ and $|T_2|$ indicate the volume of the tetrahedra $T_1$ and $T_2$, respectively, and $|f_{1,2}|$ the surface area of the face $f_{1,2}$. For $\x \in \partial T_1 \cup \partial T_2$, we have $\langle \mathbf{w}_k ( \x ), \mathbf{n}_{f_{1,2}} \rangle = \mathbf{1}_{f_{1,2}} ( \x )$.
\begin{figure}[ht]
 \begin{center}
%  \begin{tabular}{m{3.5cm}m{3.5cm}}  {\includegraphics[width=0.2\textwidth]{./img/whitney-dipole}} & 
   {\includegraphics[width=0.2\textwidth]{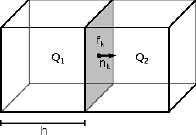}}
%  \end{tabular}
\end{center}
\caption{Zeroth-order Raviart-Thomas basis function supported 
%on two tetrahedra $T_1$ and $T_2$ (left) and 
on two hexahedra $Q_1$ and $Q_2$.}
\label{fig:rt0-support}
\end{figure}
For a regular, hexahedral mesh with edge length $h$, a $RT_0$ basis function $\mathbf{w}_k$ is supported on the two hexahedra ${Q}_1, {Q}_2 \in \mathcal{T}_h$ sharing the face $f_k = \bar{Q}_1 \cap \bar{Q}_2$ with normal vector $\mathbf{n}_k$ and centroid $\bar{\x}_k$. It can be defined via

\begin{equation}
 \mathbf{w}_k(\x) = \left\{
 \begin{aligned}
  \left( 1 - \frac{ |\langle \x - \bar{\x}_k , \mathbf{n}_k \rangle|}{h} \right) \mathbf{n}_k, &\text{ if } \x \in \bar{Q}_1 \cup \bar{Q}_2,\\
  0, \, &\text{ otherwise.}
 \end{aligned}
 \right.
\end{equation}

This definition can be transferred to nondegenerated parallelepipeds using a \emph{Piola transformation} to preserve the normal components\mbox{\cite{ned1980mixed,bre2012mixed,ber2005approximation}}. We denote the set of Raviart-Thomas basis functions $\mathbf{w}_k$ by $S^{RT_0}_h$.

The discrete approximation of \eqref{eq:mixed-op} can now be written as a matrix equation:

\begin{equation}
\label{eq:mixed-disc}
\underbrace{
 \begin{pmatrix}
  A & B^T \\ B & 0
 \end{pmatrix}
 }_{= K}
\begin{pmatrix}
 j \\
 u
\end{pmatrix}
=
\begin{pmatrix}
 b \\
 0
\end{pmatrix}
\end{equation}

with

\begin{align}
 A_{i,j} &= \int_\Omega \langle \sigma^{-1} \mathbf{w}_i , \mathbf{w}_j \rangle \dx \qquad
 B_{k,j} = \int_\Omega v_k (\nabla \cdot \mathbf{w}_j) \dx \label{eq:mixed-b}&\\
 b_i &= \int_\Omega \langle \sigma^{-1} \j^p , \mathbf{w}_i \rangle \dx & \label{eq:mixed-rhs}\\
 \intertext{\raggedleft for $v_k \in S_h^{P_0}, \, \mathbf{w}_i,\mathbf{w}_j \in S_h^{RT_0}$. $\qquad$} \nonumber
\end{align}

For the submatrices $A$ and $B$, we have $m_A = n_A = \text{\#faces}$ and $m_B = \text{\#elements}$, $n_B = \text{\#faces}$, respectively, and thus the dimension of $K$ is $m_K = n_K = \text{\#faces} + \text{\#elements}$. Using $S_h^{RT_0}$ for the matrix setup in \eqref{eq:mixed-disc}, we did not enforce the Neumann boundary condition \eqref{eq:forward-system3} in the discrete equation system so far. This has to be done explicitly when solving \eqref{eq:mixed-disc} by eliminating the respective degrees of freedom.

\subsection{Comparison to Other FE Methods for Solving the EEG Forward Problem}
The state-of-the-art FE method to solve the EEG forward
problem is the CG-FEM, for which a variety of different source models has been derived\mbox{\cite{CHW:Buc97,CHW:Yan91,CHW:Pur2011,CHW:Wol2007e,vor2016,CHW:Wol2007b}}. In addition, in\mbox{\cite{eng2015subtraction}} a discontinuous Galerkin (DG)-FEM for the EEG forward
problem has been proposed. The DG-FEM, like the Mixed-FEM, is current preserving and was derived to prevent skull leakages and to obtain more accurate and reliable results. However, whereas the Mixed-FEM actually preserves the physical current, $\mathbf{j}_h = \sigma \nabla u_h$, the DG-FEM preserves $\mathbf{j}_h=\av{\sigma \nabla u_h} - \eta \frac{\hat{\sigma}_\gamma}{\hat{h}_\gamma} \JP{u_h} \mathbf{n}$ at each element boundary, which converges to the physical current for $h\rightarrow 0$. Here, $\av{\cdot}$ and $\JP{\cdot}$ indicate the average and jump of the limit values from both sides at an (element) boundary $\gamma$, $\eta$ is a regularization parameter, and $\hat{\sigma}_\gamma$ and $\hat{h}_\gamma$ are local definitions of electric conductivity and mesh width at the surface $\gamma$\mbox{\cite{AErn2,houston2012anisotropic,eng2015subtraction,vor2016}}.

For sufficiently regular solutions, all three methods are consistent
with the strong problem and show optimal convergence rates, i.e., $O(h^2)$
in the $L^2$-norm and $O(h)$ in the energy norm for CG- and DG-FEM and $O(h)$ in the $L^2$-norm for the Mixed-FEM. Furthermore, the Mixed-FEM
and the DG-FEM are locally charge preserving. For details, we refer the reader to\mbox{\cite{CHW:Bra2007}} for the CG-FEM, to\mbox{\cite{ern2008discontinuous}} for the DG-FEM, and to\mbox{\cite{bre2012mixed}} for the Mixed-FEM.

\begin{rema}
The above-mentioned a priori convergence results will in general not apply in our case, as the dipole on the right-hand side is not in $L^2(\Omega)$. For classical, global convergence results for the CG-FEM and singular right-hand sides, we refer the reader to\mbox{\cite{cas1985singular,sco1973singular}}.
\end{rema}

CG- and DG-FEM will be used to evaluate the numerical accuracy of the approaches based on the Mixed-FEM in the numerical evaluations in Section \ref{sec:mfem-evaluation}.

\subsection{Solving the Linear Equation System \eqref{eq:mixed-disc}}
\label{sub:solving}
Due to the size of the matrix $K$ in \eqref{eq:mixed-disc}, the application of direct solvers is not feasible. Since the matrix $K$
% is indefinite as no ground potential was fixed until now and
 has a large $0$-block, Krylov subspace algorithms, such as variants of the conjugate gradient (CG) or generalized minimal residual (GMRES) method, are also not as efficient as for many other problems, since the commonly used methods for preconditioning fail\mbox{\cite{ber1994mixed}}. Nevertheless, much research has been performed to find preconditioning techniques that enable a solution using CG-solvers \cite{axe1994iterative,gol1989matrix}. A further approach to solve \eqref{eq:mixed-disc} was proposed based on the idea of introducing Lagrangian multipliers to achieve the inter\deleted{-}element continuity of the $RT_0$-basis functions, instead of including this condition by construction \cite{de1968equilibrium}. \replaced{This approach}{It} has the advantage that the resulting equation system has only \#faces unknowns, but the derivation is rather technical \cite{de1968equilibrium,bre2012mixed}. For our first evaluation of Mixed-FEM to solve the EEG forward problem, we therefore chose to apply a more direct approach that makes use of the fact that $A$ is -- unlike $K$ -- positive (semi-) definite. The chosen approach follows the ideas of\mbox{\cite{ber1994mixed}} and is based on a modification of the frequently described \emph{Uzawa-iteration}\mbox{\cite{arr1972studies,CHW:Bra2007}}. It was shown that this approach is competitive with regard to computation time when compared to the approach based on Lagrangian multipliers, called mixed-hybrid formulation in \cite{ber1994mixed}, and a (preconditioned) Augmented Lagrangian approach \cite{glo1989augmented,ber1994mixed}, in a similar scenario as the one considered here.
The origin of the derivation is identical to that of the Uzawa-iteration:
%todo cw: formulierung

 If we write \eqref{eq:mixed-disc} as a system of two equations,
\begin{subequations}
\label{eq:mixed-disc-system}
\begin{align}
 Aj + B^Tu &= b \label{eq:mixed-disc-1}\\
 Bj &= 0, \label{eq:mixed-disc-2}
\end{align}
\end{subequations}
we can left-multiply $A^{-1}$ to \eqref{eq:mixed-disc-1} and solve for $j$, i.e., $j = A^{-1}(b - B^{T} u)$. Substituting this representation of $j$ into \eqref{eq:mixed-disc-2} leads to
\begin{equation}
\begin{aligned}
 \label{eq:schur}
 Bj = BA^{-1}(b - B^Tu) &= 0 \\
\Leftrightarrow BA^{-1}B^Tu &= BA^{-1} b.
\end{aligned}
\end{equation}
$S = BA^{-1}B^T$ is the so-called \emph{Schur complement}, $m_S = n_s = \text{\#elements}$. 
 $S$ is positive semidefinite (if $\ker(B)=\{0\}$ positive definite) and since $A$ is symmetric, also $S$ is symmetric\mbox{\cite{bre2012mixed}}.
Thus, with $h = BA^{-1} b$, solving \eqref{eq:mixed-disc} is reduced to solving
\begin{equation}
\label{eq:mixed-schur}
 S u = h.
\end{equation}

\eqref{eq:mixed-schur} could now be solved using the (conjugated) Uzawa-iteration \cite{arr1972studies,bre2012mixed,CHW:Bra2007}.

However, $A^{-1}$ is a dense matrix, so that an explicit computation of $A^{-1}$ (and $S$) is not efficient considering the matrix dimensions occurring in our scenario. Instead, we access $A^{-1}$ on-the-fly by solving an additional linear equation system for each iteration, i.e., instead of calculating $x = A^{-1} y$ we solve $Ax = y$. This equation system can, for example, be solved efficiently using preconditioned CG-solvers. With the obtained implicit representation of $S$, common solver schemes such as the gradient descent or CG method can be applied to \eqref{eq:mixed-schur}.
 
When solving \eqref{eq:mixed-schur} via the CG algorithm with the implicit representation of $S$, preconditioning is advisable, as $S$ has a large condition number\mbox{\cite{ber1994mixed}}. Since $S$ is not directly accessible, it is necessary to use an approximation of $S$ for preconditioning.
%A frequent proposal for this kind of saddle point problems is to use the mass matrix of the potential space (also called pressure space due to the application in flow dynamics), i.e., $M_{i,j} = \int_\Omega \sigma v_i v_j \dx, \; v_i,v_j \in P_0(\mathcal{T}_h)$, as preconditioner \cite{elm1994inexact}. However, due to the choice of $P_0$ as scalar test-space this matrix is diagonal in our case and not expected to be an efficient preconditioner as it does not approximate $S$ well.
The use of $BB^T$ is proposed in \cite{ng1993solution}, but is efficient only in the case of constant conductivities \cite{ber1994mixed}. Although $BB^T$ approximates the pattern of $S$ well enough, it does not provide a reasonable approximation of the matrix entries of $S$. Instead, it is suggested to choose a diagonal matrix $D$ that in some sense approximates $A$ and to use $BD^{-1}B^T$ as input to the preconditioner \cite{ber1994mixed}. It is further proposed to choose $D_{i,i} = l_2(A_{i,:}) = (\sum_j A_{i,j}^2 )^{1/2}$.
Indeed, this approximation led to the best results when it was compared to the choices $D_{i,i} = A_{i,i}$, $D_{i,i} = \sum_j A_{i,j}$, and $D_{i,i} = l_1(A_{i,:}) = \sum_j |A_{i,j}|$ \cite{vor2016}. % tough the differences in performance between defining $D_{i,i}$ as the $l_1$- or $l_2$-column norm were relatively small.

Since all considered choices for $D$ are diagonal, the structure of the matrix $P = BD^{-1}B^T$ is identical to the structure of $BB^T$ and cannot be easily inverted. Also, due to the structure of $P$, commonly chosen preconditioners such as the incomplete LU-factorization (ILU) cannot be expected to be efficient \cite[p. 330]{che2005matrix}. We found that approximating $P^{-1}$  
% by an iterative solution using an ILU(0)-CG solver and
using an algebraic multigrid (AMG) method leads to a  performance that is sufficient for our first evaluations\cite{vor2016}.
%, with the AMG-preconditioner being fastest.

Besides preconditioning of the ``outer iteration'', a further speed-up of the solver could be achieved by reducing the accuracy with which the inner equation, $Ax = y$, is solved. This approach can be interpreted to be similar to inexact Uzawa-algorithms as they are proposed in the literature \cite{elm1994inexact}. Since reducing the number of iterations for solving the inner equation did not result in an increase in the number of outer iterations that is necessary to reach the desired solution accuracy, performing only one iteration led to the fastest solving speed. Using this approach, solving the equation system \eqref{eq:mixed-schur} took less than two minutes for the finest used spherical model with 1 mm mesh resolution (model \segres{1}{1} in Table \ref{tab:models}) \cite{vor2016}.

Through the integration of algebraic multigrid preconditioners to the Uzawa-like method proposed in \cite{ber1994mixed}, our solution algorithm has similarities to the combined conjugate gradient-multigrid algorithm proposed in \cite{ver1984combined}. However, in \cite{ver1984combined} no preconditioning of the outer iteration is performed.

\subsection{Modeling of a Dipole Source}
This section focuses on the exact choice of the source distribution $\j^p$. In principle, arbitrary distributions $\j^p \in \L2^3$, $\supp \j^p \subset \Omega^\circ$, can be modeled. The common choice in EEG forward modeling is $\j^p = \mathbf{m} \delta_{\x_0}$, where $\delta_{\x_0}$ is the Dirac delta distribution and $\mathbf{m}$ the dipole moment. Since maximally $\delta \in H^{- 3/2 - \epsilon}$, the assumption $\j^p \in \L2^3$ is violated. The authors are not aware of any literature investigating the influence of singular right-hand sides $\j^p$ for the Mixed-FEM. However, in the case of the CG-FEM, it was shown that such a singular right-hand side does not affect the existence and uniqueness of a solution in general, but leads to a lower regularity of the solution, and, in consequence, to worse global a priori error estimates\mbox{\cite{cas1985singular,sco1973singular}}. (Quasi-) optimal convergence for the CG-FEM can be shown in seminorms that exclude the locations of the singularities\mbox{\cite{kop2014optimal}}.

As \eqref{eq:forward} is represented by a system of two PDEs now, there are two options to model the dipole source. The dipole can be modeled either in the ``current space'' \eqref{eq:mixed-weak1} or in the ``potential space'' \eqref{eq:mixed-weak2} (sometimes also called ``pressure space'' due to the origin of Mixed-FEM in reservoir simulations \cite{ng1993solution}). The first option corresponds to an evaluation of the functional $l$ in the discrete space $RT_0$ as it was defined in \eqref{eq:mixed-rhs}. For $\j^p = \mathbf{m} \delta_{\x_0}$, i.e., a current dipole with moment $\mathbf{m}$ at position $\x_0$, we have
\begin{align}
\label{eq:rt-b-direct}
  b_i = b^{cur}_i &= \int_\Omega \langle \sigma^{-1} \mathbf{m}\delta_{\x_0}, \mathbf{w}_i \rangle \dx \nonumber\\
  &= \left\{ 
  \begin{aligned}
   \langle \sigma^{-1} \mathbf{m} , \mathbf{w}_i (\x_0) \rangle, &\text{ if } \x_0 \in \supp \mathbf{w}_i \\
   0, &\text{ otherwise.}
  \end{aligned} \right.
\end{align}
This approach will be called the \emph{direct approach} with $h = h^{direct} = BA^{-1} b^{cur}$.

\begin{figure}[htb]
 \begin{center}
  {\includegraphics[width=0.441\textwidth]{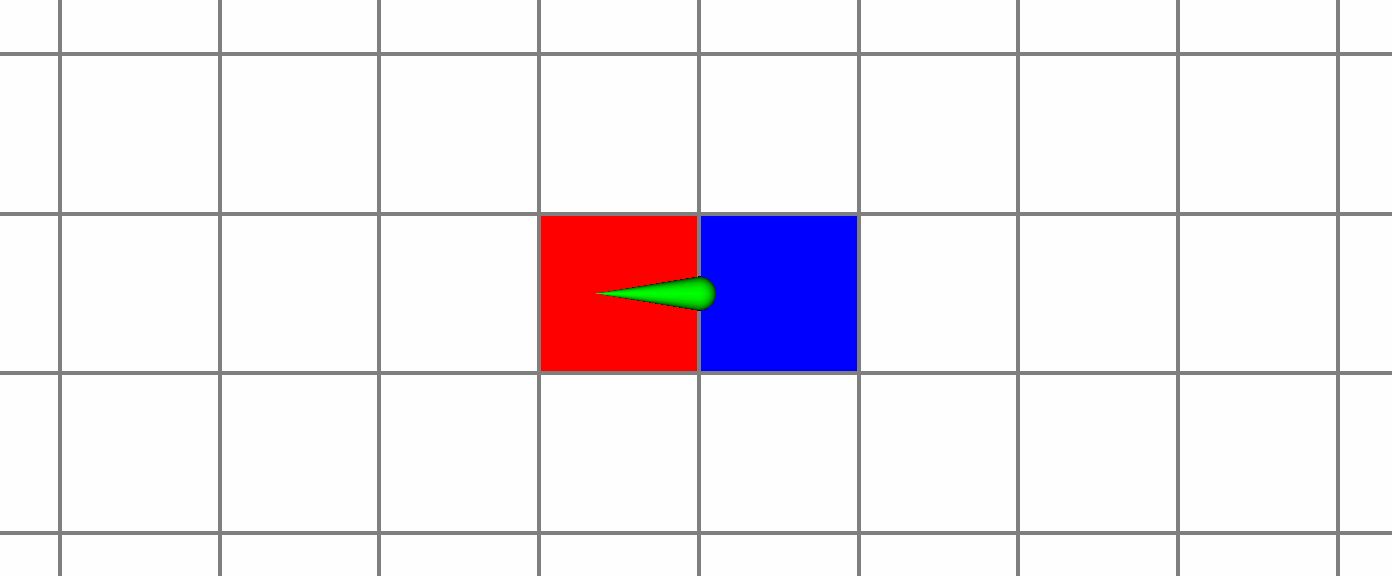}} \\
  \vspace{.01\textwidth}
    {\includegraphics[width=0.441\textwidth]{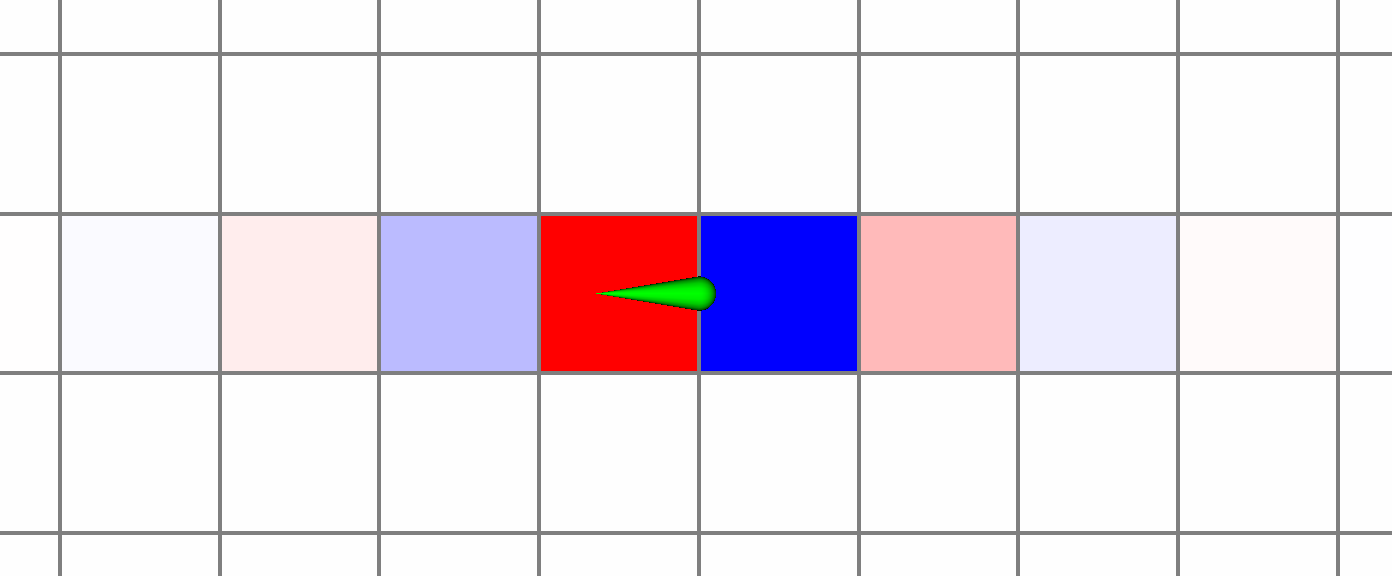}} \\
      \vspace{.01\textwidth}
       {\includegraphics[width=0.441\textwidth]{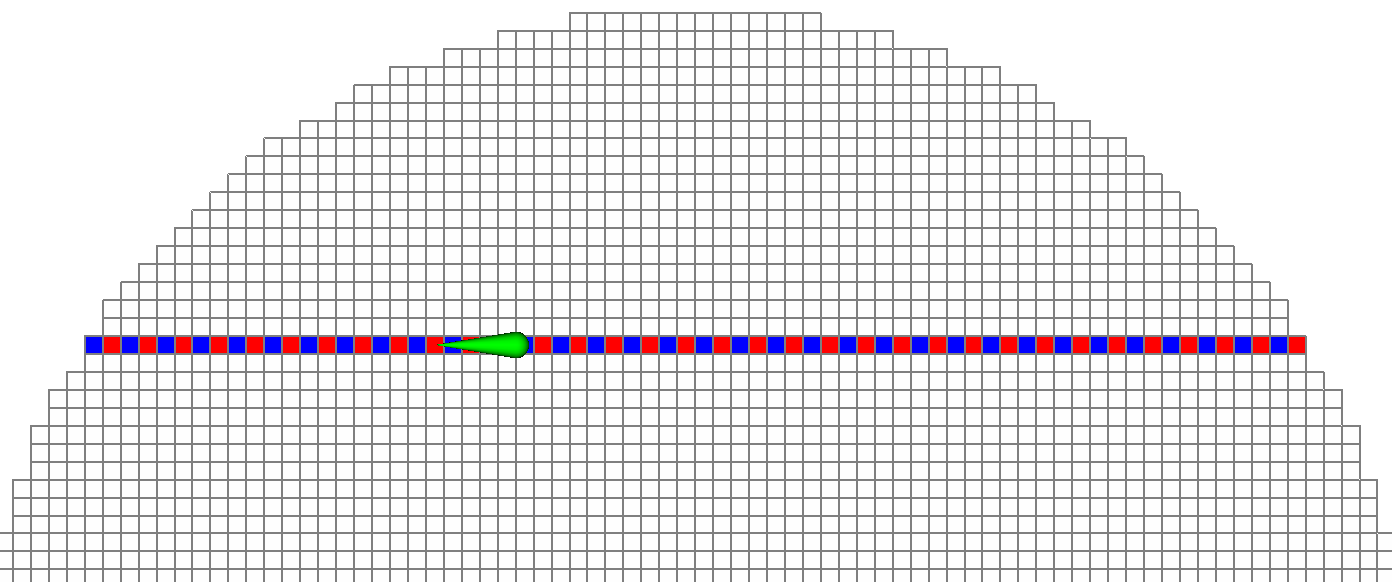}}
  \end{center}
\caption{Visualization of $h = h^{proj} = B b^{pot}$ (top), $h = h^{direct} = BA^{-1} b^{cur}$ (middle), and full view of the support of $h^{direct}$ through visualizing $\sign(h^{direct})$ (bottom) for a source positioned in the center of a face $f_i$ and direction $\mathbf{n}_{f_i}$ (green cone). The slice is taken at the dipole position in the y-plane. The coloring indicates the values for the $P_0$ basis function corresponding to the respective element; red is positive, blue is negative.}
\label{fig:mixed-rhs-fitted}
\end{figure}

A representation of the dipole in the potential space, henceforth called the \emph{projected approach}, can be obtained using the matrix $B$, which can be interpreted as a mapping between the current and the potential space. Figuratively, the (source) current is mapped to the distribution of sinks and sources generating this current. The projected approach is similar to the Whitney approach that was introduced for the CG-FEM \cite{CHW:Pur2011,bau2015comparison}, except for using the scalar space $P_0$ instead of $P_1$. In both approaches, a current source, represented by $RT_0$ basis functions, is mapped to the potential space. To achieve this representation for the Mixed-FEM, we redefine $b$ to be the approximation of $\j^p$ in the space $RT_0$

\begin{align}
\label{eq:rt-b-direct-pot}
  b^{pot}_i &= \int_\Omega \langle \mathbf{m}\delta_{\x_0}, \mathbf{w}_i \rangle \dx \nonumber\\
  &= \left\{ 
  \begin{aligned}
   \langle \mathbf{m} , \mathbf{w}_i (\x_0) \rangle, &\text{ if } \x_0 \in \supp \mathbf{w}_i \\
   0, &\text{ otherwise.}
  \end{aligned} \right.
\end{align}
 $b^{pot}$ is then projected to the space $P_0$ using $B$. We obtain $h = h^{proj} = B b^{pot}$; the dipole is represented by a source and a sink in the potential space in this case (Figure \ref{fig:mixed-rhs-fitted}, top).

\begin{rema}
If a single $RT_0$ function is chosen as the source distribution and a hexahedral mesh is used, i.e., the source is positioned on the face $f_i$ and the direction is $\mathbf{n}_{f_i}$, only one entry of $b$ is nonzero (cf. \eqref{eq:rt-b-direct}, \eqref{eq:rt-b-direct-pot}). When applying the projection to the potential space using the matrix $B$, which has only two nonzero entries per column (cf. \eqref{eq:mixed-b}), the right-hand side vector, which is given by $h = h^{proj} = B b^{pot}$, also has only two nonzero entries (Figure \ref{fig:mixed-rhs-fitted}, top). In contrast, the right-hand side $h^{direct} = BA^{-1} b^{cur}$ causes a blurring of the current source when interpreting it as a monopole distribution and visualizing it in the pressure space. It leads to nonzero right-hand side entries $h_i$ assigned to all elements that are ``in the source direction'' (cf. Figure \ref{fig:mixed-rhs-fitted}, middle and bottom; Figure \ref{fig:mixed-rhs-fitted}, bottom shows the sign function of all elements corresponding to non\deleted{-}zero right-hand side entries through red-blue coloring). However, most of these values are small.
% as the truncated visualization for arbitrary dipole position and direction indicates (Figure \ref{fig:mixed-rhs}, middle).
% (Figure \ref{fig:mixed-rhs}).

This structure of $b$ transforms accordingly to the case of arbitrarily positioned and oriented sources, as the right-hand side vectors $b$ -- and thereby also $h$ -- are linear combinations of the solutions for dipoles oriented in the directions of the mesh basis vectors in this case. The accuracies of the different representations are evaluated in Section \ref{sec:mfem-evaluation}.
\end{rema}

%\begin{figure*}[htb]
% \begin{center}
%  {\includegraphics[width=0.29\textwidth]{./img/mfem_projected}} \hfil
%    {\includegraphics[width=0.29\textwidth]{./img/mfem_direct_trunc_5}} \hfil {\includegraphics[width=0.29\textwidth]{./img/mfem_direct_full3}}
%  \end{center}
%\caption[M-FEM right-hand sides for random dipole]{Visualization of $h= h^{proj}$ (left), $h^{direct} = BA^{-1} b^{direct}$ truncated at 5\% of maximum (middle), and full visualization of $h^{direct} = BA^{-1} b^{direct}$ (right) for dipole with random position and orientation. The coloring indicates the values for the $P_0$ basis function corresponding to the respective element; red is positive, blue is negative.}
%\label{fig:mixed-rhs}
%\end{figure*}

\section{Methods}
\subsection{Implementation}
\IEEEPARstart{F}{or} this study, both the direct (i.e., $h = h^{direct} = BA^{-1} b^{cur}$) and the projected ($h = h^{proj} = B b^{pot}$) Mixed-FEM approaches were implemented in the
DUNE framework \cite{dune08:1,dune08:2} using the DUNE-PDELab toolbox
\cite{bastian2010generic}. In addition, a solver corresponding to a conjugate Uzawa-iteration with additional preconditioning and implicit representation of $A^{-1}$, as derived in Section \ref{sub:solving}, was implemented using the CG-solver template from the DUNE module \emph{iterative solvers template library} (DUNE-ISTL) in combination with %different preconditioners, e.g., AMG, ILU(0), and iterative.
the AMG preconditioner \cite{bla2010parallel}.

\subsection{Evaluation}
In order to evaluate the accuracy of the Mixed-FEM, different comparisons both in hexahedral four-layer sphere models and in realistic head models were performed. As is common for the evaluation of EEG forward approaches, the error measures RDM (minimal error 0, maximal error 2) and lnMAG (minimal error 0, maximal error $\pm \infty$)  were used \cite{CHW:Mei89,CHW:Gue2010}.
\begin{equation}
\label{eq:rdm-lnmag}
\begin{split}
  RDM ( u^{num}, u^{ref} ) &= \left\| \frac{u^{num}}{\| u^{num} \|_2} - \frac{u^{ref}}{\| u^{ref} \|_2}  \right\|_2 \\
 lnMAG ( u^{num}, u^{ref} ) &= \ln \left(\frac{\| u^{num} \|_2}{\| u^{ref} \|_2}\right)
 \end{split}
\end{equation}
 In the sphere models, the solution was evaluated on the whole outer boundary instead of using single electrode positions, so that the results are independent of the choice of sensor positions. For the realistic head model, the sensor positions of a realistic 80-electrode EEG cap were used \cite{CHW:Vor2014,vor2016}.
 
\begin{table}[htbp]
\renewcommand{\arraystretch}{1.3}
\caption[Four-layer Sphere Models]{Four-layer sphere models (compartments from in- to outside)}
\label{tab:4-layer-compartments}
\centering
\begin{tabular}{lrrr}
\hline
Compartment & \multicolumn{1}{l}{Outer Radius} & \multicolumn{1}{l}{$\sigma$} & \multicolumn{1}{l}{Reference}\\
\hline \hline
Brain & 78 mm & 0.33 S/m & \cite{CHW:Ram2004} \\
CSF & 80 mm & 1.79 S/m & \cite{CHW:Bau97} \\
Skull & 86 mm & 0.01 S/m & \cite{CHW:Dan2011} \\
Skin & 92 mm & 0.43 S/m & \cite{CHW:Dan2011,CHW:Ram2004} \\
\hline
\end{tabular}
\end{table}

\begin{table}[tb]
\renewcommand{\arraystretch}{1.3}
\caption[Sphere Model Parameters]{Sphere Model Parameters}
\label{tab:models}
\centering
\begin{tabular}{lrrrr}
\hline
& \multicolumn{1}{l}{Mesh width ($h$)} & \multicolumn{1}{l}{\#vertices} & \multicolumn{1}{l}{\#elements} & \multicolumn{1}{l}{\#faces}\\
\hline \hline
\segres{1}{1} & 1 mm & 3,342,701 & 3,262,312 & 9,866,772 \\
\segres{2}{2} & 2 mm & 428,185 & 407,907 & 1,243,716\\
\hline
\end{tabular}
\end{table}

\begin{table}[tb]
\renewcommand{\arraystretch}{1.3}
\caption[Realistic Head Model Parameters]{Realistic Head Model Parameters}
\label{tab:models-realistic}
\centering
\begin{tabular}{lrrrr}
\hline
& \multicolumn{1}{l}{Mesh width ($h$)} & \multicolumn{1}{l}{\#vertices} & \multicolumn{1}{l}{\#elements} & \multicolumn{1}{l}{\#faces}\\
\hline \hline
\ci{6}{hex}{1mm} & 1 mm & 3,965,968 & 3,871,029 & 11,707,401\\
\ci{6}{hex}{2mm} & 2 mm & 508,412 & 484,532 & 1,477,164\\
\ci{6}{tet}{hr} & -- & 2,242,186 & 14,223,508 & 27,314,610\\
\hline
\end{tabular}
\end{table}

\begin{table}[b]
\caption{Model leaks}
\label{tab:leaks}
\centering
\begin{tabular}{lcrrr}
\hline
Model & \multicolumn{1}{l}{Outer Skull Radius} & \multicolumn{1}{l}{\#leaks} \\
\hline \hline
\segresR{2}{2}{82} & 82 mm & 10,080 \\
\segresR{2}{2}{83} & 83 mm & 1,344 \\
\segresR{2}{2}{84} & 84 mm & 0 \\
\hline
\end{tabular}
\end{table}

Besides the two Mixed-FEM approaches, the Whitney CG-FEM was included in our sphere model comparisons, as it relies on the same approximation of the dipole source \cite{CHW:Pur2011,bau2015comparison}. By including the Whitney CG-FEM, the differences between Mixed- and CG-FEM can be directly evaluated. Two four-layer hexahedral sphere models, \segres{1}{1} and \segres{2}{2}, with a mesh resolution of 1 and 2 mm, respectively, were generated (Tables \ref{tab:4-layer-compartments}, \ref{tab:models}). Sources were placed at 10 different radii, and for each radius 10 sources were randomly distributed. This distribution of the test sources allows us to gain a statistical overview of the range of the numerical accuracy at each eccentricity. Since the numerical errors increase along with the eccentricity, i.e., the quotient of source radius and radius of the innermost compartment boundary, the radii of the source positions were chosen so that the distances to the next conductivity jump (brain/CSF boundary) were logarithmically distributed. The most exterior eccentricity 0.993 corresponds to a distance of only $\approx$ 0.5 mm to the conductivity jump. In praxis (and for the realistic head model used in this study), sources are usually placed so that at least one layer of elements is between the source element and the conductivity jump, which is fulfilled for the considered eccentricities $\leq 0.987$ in the 1 mm model and the eccentricities $\leq 0.964$ in the 2 mm model. The reference solutions $u^{ref}$ were computed using a quasianalytical solution for sphere models \cite{mun1993}.

In the first study, for each model, the sources were placed on the closest face center and the source directions were chosen according to the face normals, so that only one basis function contributes to the right-hand side vectors $b$ (cf. \deleted{Equations }\eqref{eq:mixed-b}, \eqref{eq:rt-b-direct}). Therefore, the results are not influenced by the interpolation that is needed for arbitrary source directions and positions. For the Whitney approach, it was shown that it has the highest accuracy of all CG-FEM approaches in this scenario \cite{bau2015comparison}. Next, the three approaches were compared in the same models using the initially generated random source positions and radial source directions, so that neither positions nor directions were adjusted to the mesh. We limit our investigations to radial sources, as eccentric radial sources were shown to lead to higher numerical errors than tangential sources in previous studies \cite{CHW:Wol2007a}. Finally, the projected Mixed-FEM and Whitney CG-FEM were evaluated in combination with the models \segresR{2}{2}{82}, \segresR{2}{2}{83}, and \segresR{2}{2}{84} generated from model \segres{2}{2} but with an especially thin skull layer, again with random positions and radial source directions. Table \ref{tab:leaks} indicates the outer skull radii of the different models and the resulting number of leakages, i.e., the number of nodes in which elements of skin and CSF compartment touch.

\begin{figure}[tb]%
\centering
\includegraphics[width=.19\textwidth]{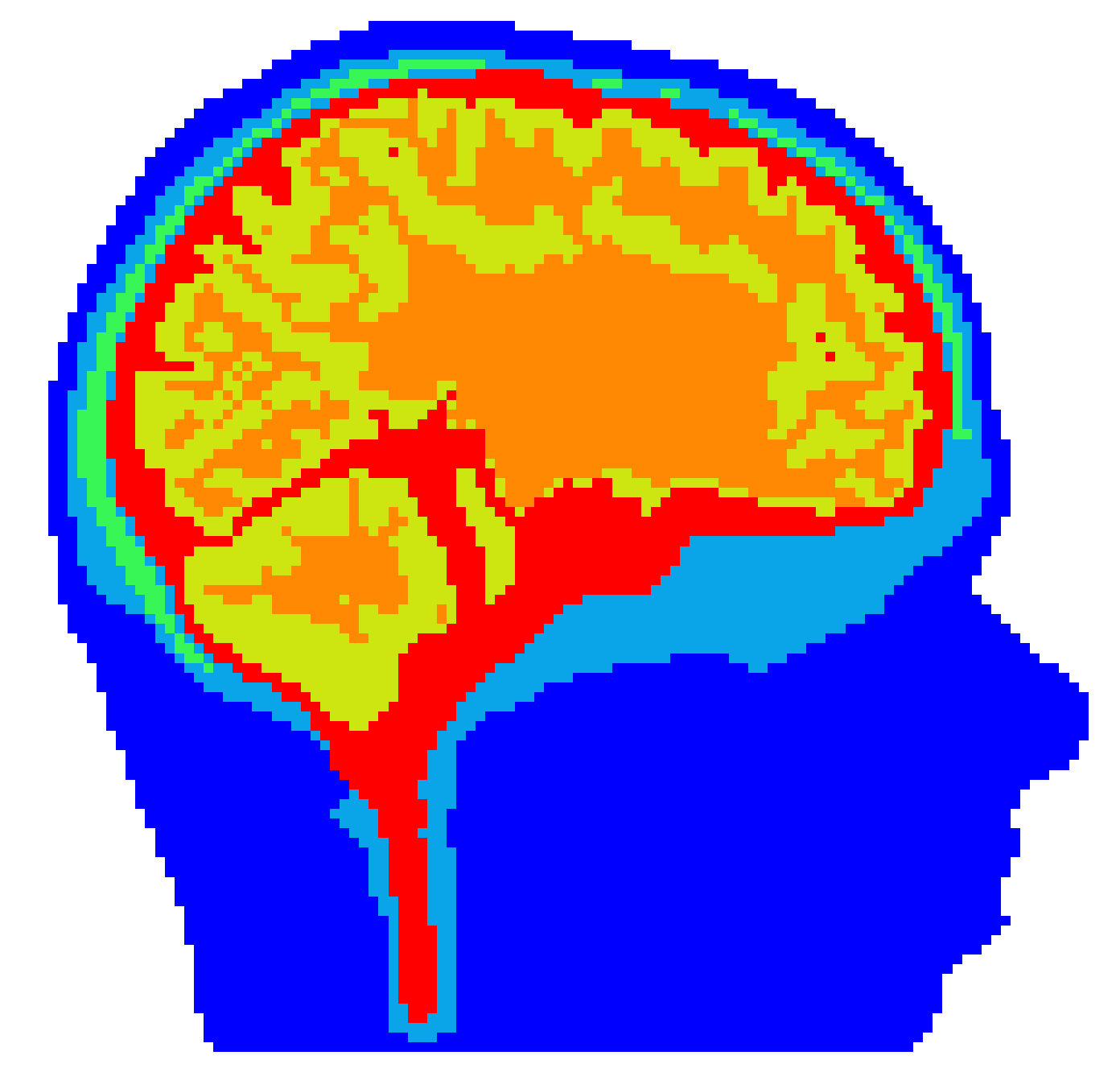} \hfil
\includegraphics[width=.19\textwidth]{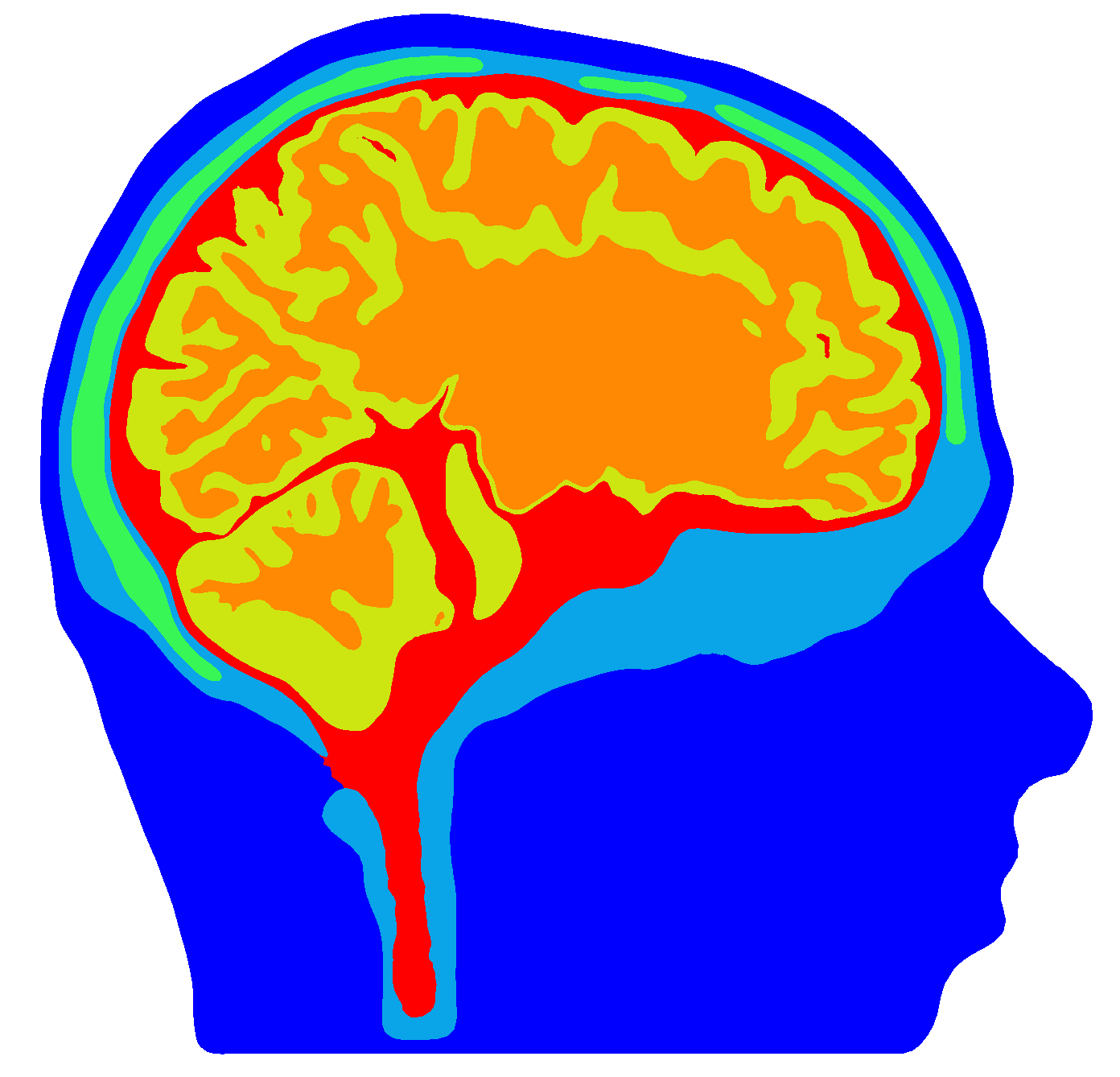}
\caption{Visualization of realistic six-compartment hexahedral (\ci{6}{hex}{2mm}, left) and high-resolution reference head model (\ci{6}{tet}{hr}, right).}%
\label{fig:headmodel}%
\end{figure}

Mixed-FEM, CG-FEM, and DG-FEM were further evaluated in a more realistic scenario. Two realistic six-compartment hexahedral head models with mesh widths of 1 mm, \ci{6}{hex}{1mm}, and 2 mm, \ci{6}{hex}{2mm}, were created, resulting in 3,965,968 vertices and 3,871,029 elements and 508,412 vertices and 484,532 elements, respectively (Table \ref{tab:models-realistic}, Figure \ref{fig:headmodel}). As the model with a mesh width of 2 mm was not corrected for leakages, 1,164 vertices belonging to both CSF and skin elements were found, mainly located at the temporal bone. The conductivities were chosen according to \cite{CHW:Vor2014}. Of 18,893 source positions placed in the gray matter with a normal constraint, those not fully contained in the gray matter compartment (i.e., where the source was placed in an element at a compartment boundary) were excluded. In consequence, 17,870 source positions remained for the 1 mm model and 17,843 source positions for the 2 mm model. As sensor configuration an 80 channel realistic EEG cap was chosen. The investigated approaches were projected Mixed-FEM, Whitney CG-FEM, St. Venant CG-FEM \cite{CHW:Buc97}, and Partial Integration DG-FEM \cite{vor2016, eng2015subtraction}. St. Venant CG-FEM  and Partial Integration DG-FEM were additionally included, since they were shown to achieve the highest accuracies of the different CG- and DG-FEM approaches, respectively, when choosing arbitrary source directions and positions \cite{eng2015subtraction,bau2015comparison}. Solutions for all methods were computed in the 2 mm model, and a solution in the 1 mm model was calculated using the St. Venant CG-FEM. In the realistic scenario, RDM and lnMAG were evaluated in comparison to a reference solution that was computed using the St. Venant method in a high-resolution tetrahedral model, \ci{6}{tet}{hr}, based on the same segmentation (Table \ref{tab:models-realistic}, 2,242,186 vertices, 14,223,508 elements). For details of this model, we refer the reader to \cite{CHW:Vor2014,vor2016}.

\section{Results}
\label{sec:mfem-evaluation}
In this paper, a new finite element method to solve the EEG forward problem is introduced. It is expected that it should be preferrable compared to the commonly used CG-FEM approaches especially in leakage and realistic scenarios. The goal of Sections \ref{sub:mfem-evaluation1} and \ref{sub:mfem-evaluation2} is to show that this new method performs appropriately when compared to the established CG-FEM in common sphere models, and in Sections \ref{sub:results-leaky} and \ref{sub:realistic} the accuracy in leakage and realistic scenarios is evaluated.

\subsection{Comparison of Whitney CG-FEM and Mixed-FEM for Optimal Source Positions}
\label{sub:mfem-evaluation1}
\begin{figure}[tb]%
\centering
\includegraphics[width=.46\textwidth]{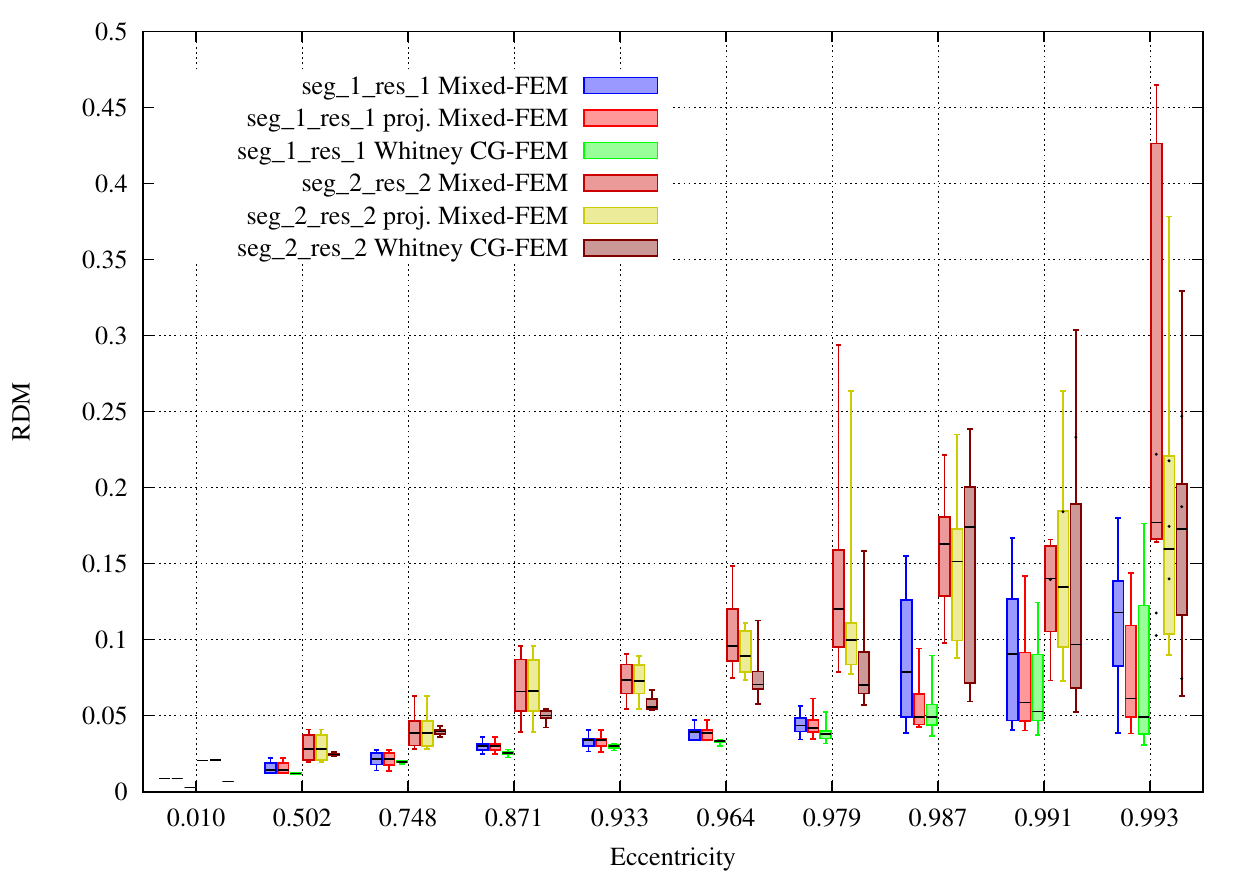}\\%
\includegraphics[width=.46\textwidth]{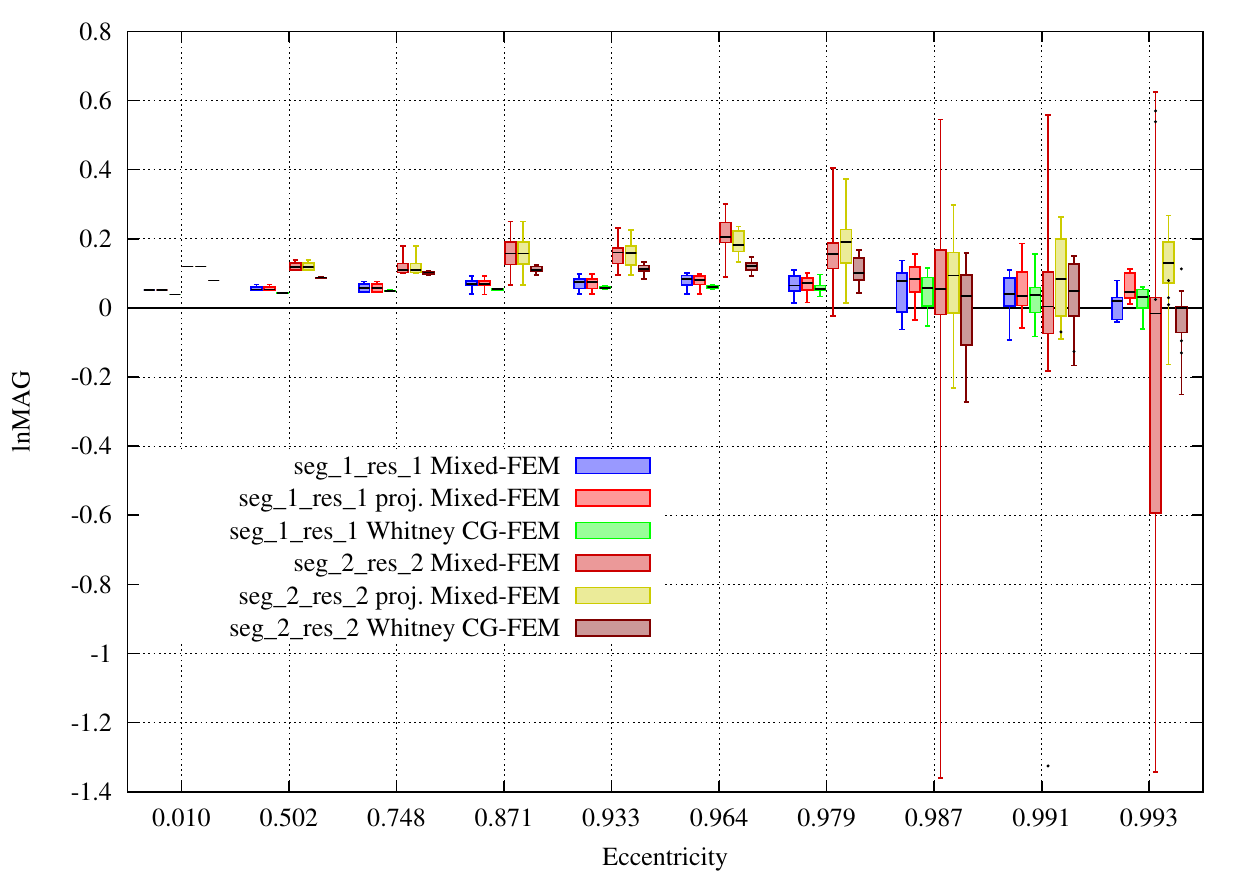}
%\caption[Comparison of direct and projected M-FEM and Whitney CG-FEM, optimal positions and directions]{Comparison of direct and projected Mixed-FEM and Whitney CG-FEM in meshes \segres{1}{1} and \segres{2}{2}. Results for optimized dipole positions. Visualized are the mean error (left column) and boxplots (right column) of RDM (top row) and lnMAG (bottom row). Dipole positions that are outside the brain compartment in the discretized models are marked as dots.
%Note the different scaling of the x-axes.}%
\caption[Comparison of direct and projected Mixed-FEM and Whitney CG-FEM, optimal positions and directions]{Comparison of direct and projected Mixed-FEM and Whitney CG-FEM in meshes \segres{1}{1} and \segres{2}{2}. Results for optimized dipole positions. Visualized boxplots of RDM (top row) and lnMAG (bottom row). Dipole positions outside the brain compartment in the discretized models are marked as dots.
Note the logarithmic scaling of the x-axes.}
\label{fig:radial-conv-mixed}%
\end{figure}
\IEEEPARstart{C}{omparing} the three approaches with regard to the RDM in model \segres{1}{1} (Figure \ref{fig:radial-conv-mixed}), no remarkable differences are found up to an eccentricity of 0.964 (distance from next conductivity jump $\geq 2.8$ mm) with maximal errors below 0.05 for all approaches (Figure \ref{fig:radial-conv-mixed}, top row). At an eccentricity of 0.979 (dist. $\approx 1.6$ mm), the maximal errors for the Mixed-FEM slightly increase. However, the maximal errors remain clearly below 0.1. Also the Whitney CG-FEM has a maximal error below 0.1 at this eccentricity, and the upper quartile and median are lower than for the Mixed-FEM. For the highest three eccentricities, the RDM clearly increases for all considered approaches. The variance, especially for the highest eccentricities, is lowest for projected Mixed-FEM and Whitney CG-FEM. In the coarser model \segres{2}{2}, direct and projected Mixed-FEM perform similar up to eccentricities of 0.933 or 0.964 (dist. $\geq 2.8$ mm), whereas the errors for the Whitney CG-FEM are lower and have less variance. For higher eccentricities, a rating of the accuracies is hardly possible due to the higher variance.

With regard to the lnMAG (Figure \ref{fig:radial-conv-mixed}, bottom row), only minor differences are recognizable for model \segres{1}{1}. In model \segres{2}{2}, it is notable that the direct Mixed-FEM leads to very high maximal errors for eccentricities of 0.987, whereas Whitney CG-FEM and projected Mixed-FEM perform similar with a tendency of the Whitney CG-FEM toward lower errors.

\subsection{Comparison of Whitney CG-FEM and Mixed-FEM for Random Source Positions}
\label{sub:mfem-evaluation2}
\begin{figure}[tb]%
\centering
 \includegraphics[width=.46\textwidth]{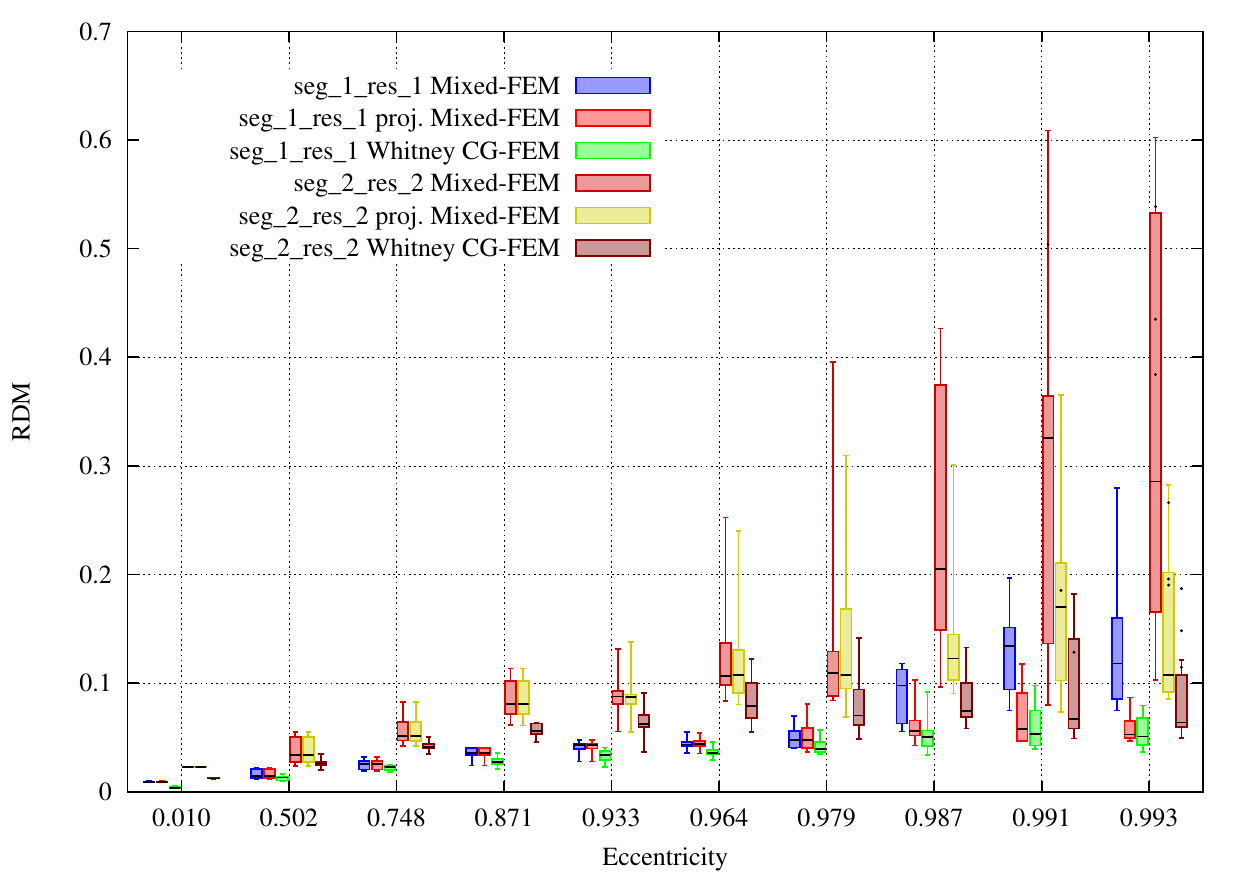}\\%
 \includegraphics[width=.46\textwidth]{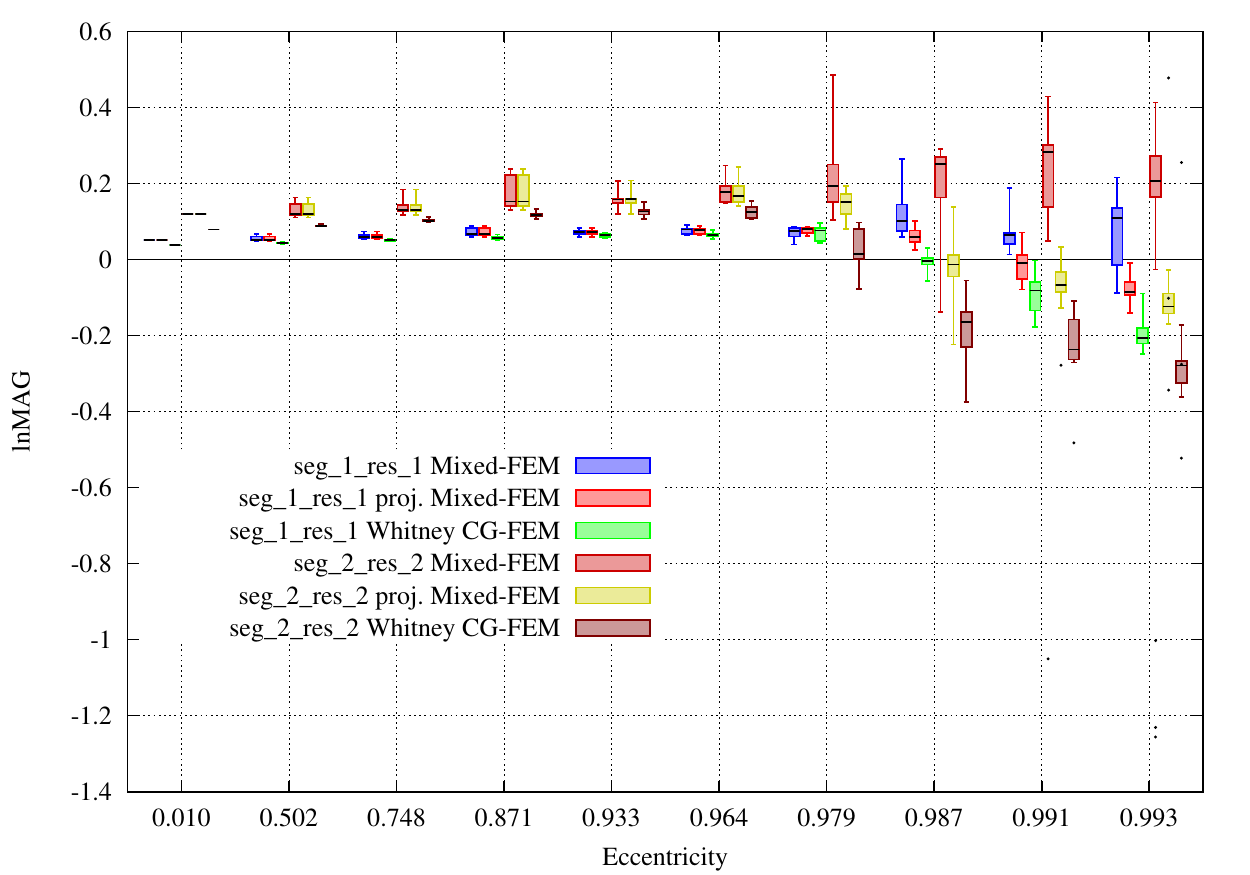}
%\caption[Comparison of direct and projected M-FEM, random positions and radial directions]{Convergence with increasing mesh and geometrical resolution for Mixed-FEM approaches. Results of radial dipole computations. Visualized are the mean error (left column) and boxplots (right column) of RDM (top row) and lnMAG (bottom row). Dipole positions that are outside the brain compartment in the discretized models are marked as dots.
%Note the different scaling of the x-axes.}%
\caption[Comparison of direct and projected Mixed-FEM and Whitney CG-FEM, optimal positions and directions]{Comparison of direct and projected Mixed-FEM and Whitney CG-FEM in meshes \segres{1}{1} and \segres{2}{2}. Results for random dipole positions. Visualized boxplots of RDM (top row) and lnMAG (bottom row). Dipole positions outside the brain compartment in the discretized models are marked as dots.
Note the logarithmic scaling of the x-axes.}
\label{fig:radial-geom-mixed}%
\end{figure}
The next comparison expands the previous results to random source positions and radial source orientations. When comparing the two Mixed-FEM approaches with regard to the RDM (Figure \ref{fig:radial-geom-mixed}, top row), both models show no major differences up to an eccentricity of 0.964 (dist. $\geq 2.8$ mm), but the Whitney CG-FEM leads to lower errors especially in model \segres{2}{2}. For model \segres{1}{1}, the RDM is constantly below 0.05 at low eccentricities (up to eccentricity $\leq 0.964$, i.e., dist. $\geq 2.8$ mm). With increasing eccentricity, the RDM for the projected Mixed-FEM and Whitney CG-FEM mainly remains below 0.1, whereas the maximal RDM is at nearly 0.3 for the direct approach and the median is above 0.1. Also in model \segres{2}{2}, the projected approach outperforms the direct approach with regard to the RDM. The less accurate approximation of the geometry leads to higher errors in these models, e.g., the minimal RDM at an eccentricity of 0.964 (dist. $\geq 2.8$ mm) is already at nearly 0.1 for both approaches in model \segres{2}{2}. The Whitney CG-FEM performs clearly better than both Mixed-FEM approaches in this model, with maximal errors below 0.13 at this eccentricity. For more eccentric sources, the projected approach, again, performs better than the direct approach. Nevertheless, the errors for the Whitney CG-FEM remain at a lower level.

The results for the lnMAG (Figure \ref{fig:radial-geom-mixed}, bottom row) do not show remarkable differences for all models up to an eccentricity of 0.964. In model \segres{1}{1}, the projected Mixed-FEM leads to the lowest spread for the three highest eccentricities. However, the lnMAG decreases from positive  values for all source positions at low eccentricities to completely negative values at the highest eccentricity. This effect is even stronger for the Whitney CG-FEM. In contrast, the median of the direct Mixed-FEM remains close to constant up to the highest eccentricity, but with a higher spread. The same behavior of the three approaches, just at a generally higher error level, is found for model \segres{2}{2}.

\subsection{Comparison of Mixed-FEM Approaches in Leaky Sphere Models}
\label{sub:results-leaky}
\begin{figure}[tb]%
\centering
 \includegraphics[width=.46\textwidth]{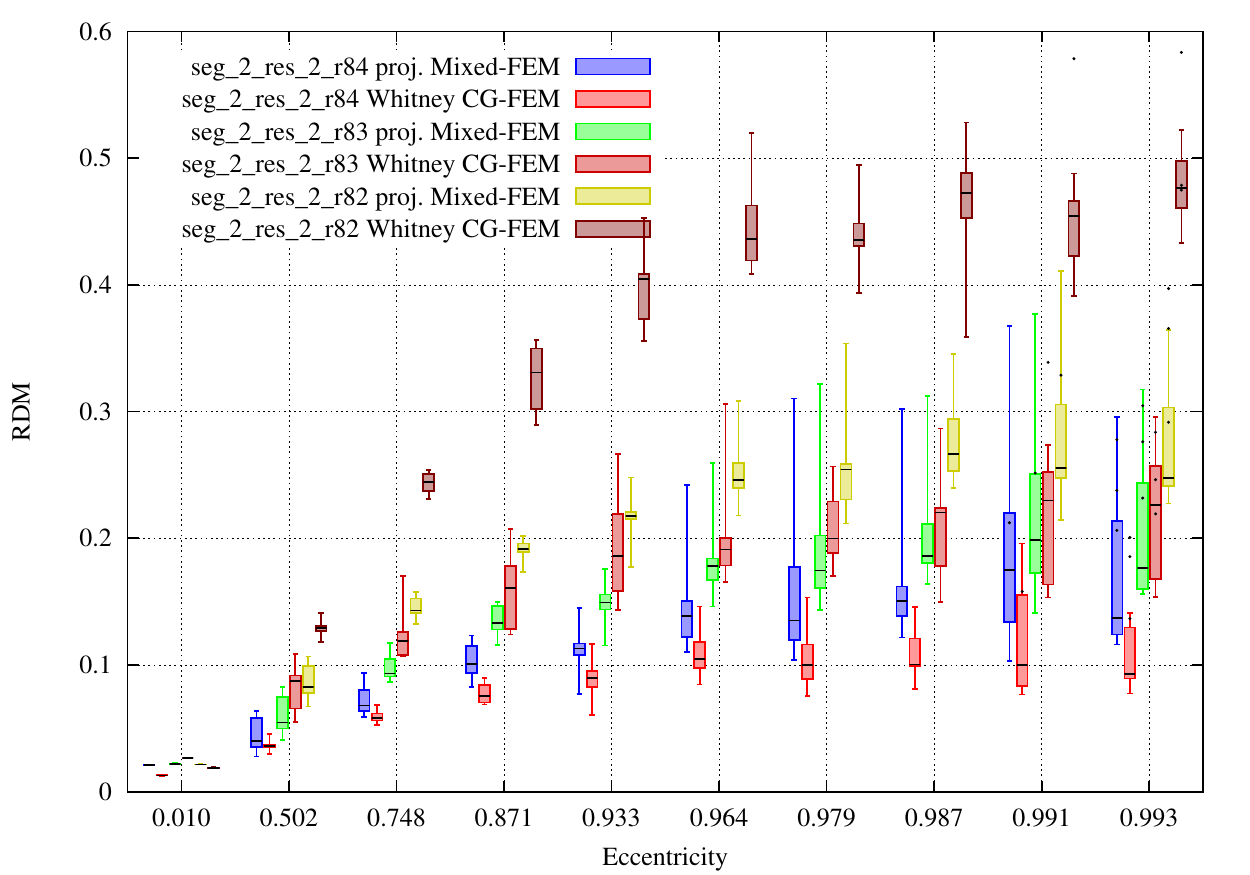}\\%
 \includegraphics[width=.46\textwidth]{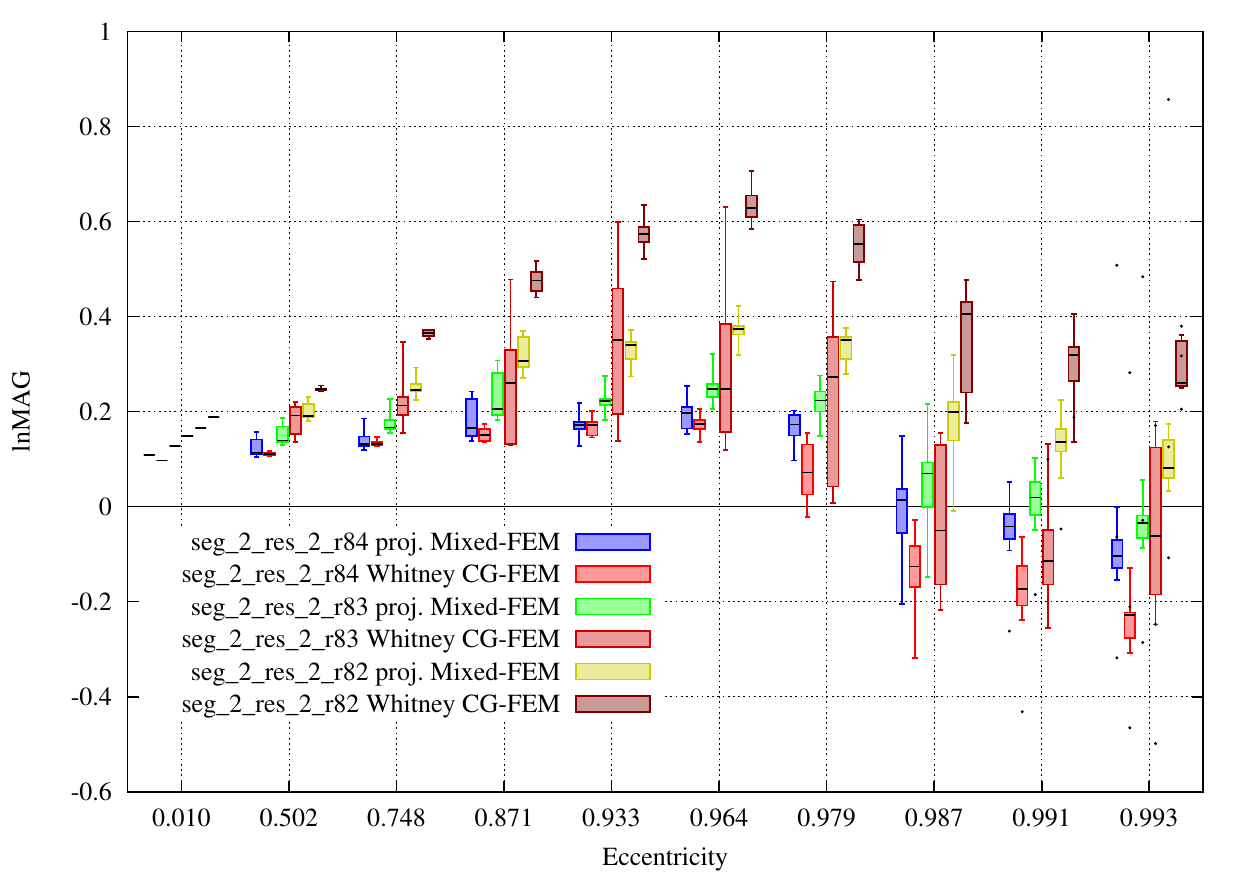}
\caption{Comparison of projected Mixed-FEM and Whitney CG-FEM in meshes with thin skull compartment. Results for random dipole positions. Visualized boxplots of RDM (top row) and lnMAG (bottom row). Dipole positions outside the brain compartment in the discretized models are marked as dots.
Note the logarithmic scaling of the x-axes.}
%[Comparison of direct and projected Mixed-FEM in models with thin skull, random positions and radial directions]{Comparison of increase of errors for decreasing skull thickness between the two Mixed-FEM approaches.
% Results of radial dipole computations. Visualized are the mean error (left column) and boxplots (right column) of RDM (top row) and lnMAG (bottom row). Dipole positions that are outside the brain compartment in the discretized models are marked as dots. Note the different scaling of the x-axes.}%
%Results of radial dipole computations. Visualized boxplots of RDM (top row) and lnMAG (bottom row). Dipole positions that are outside the brain compartment in the discretized models are marked as dots.
%Note the logarithmic scaling of the x-axes.}
\label{fig:radial-leak-mixed}%
\end{figure}
\begin{figure*}[tb]%
\centering
\includegraphics[width=.24\textwidth]{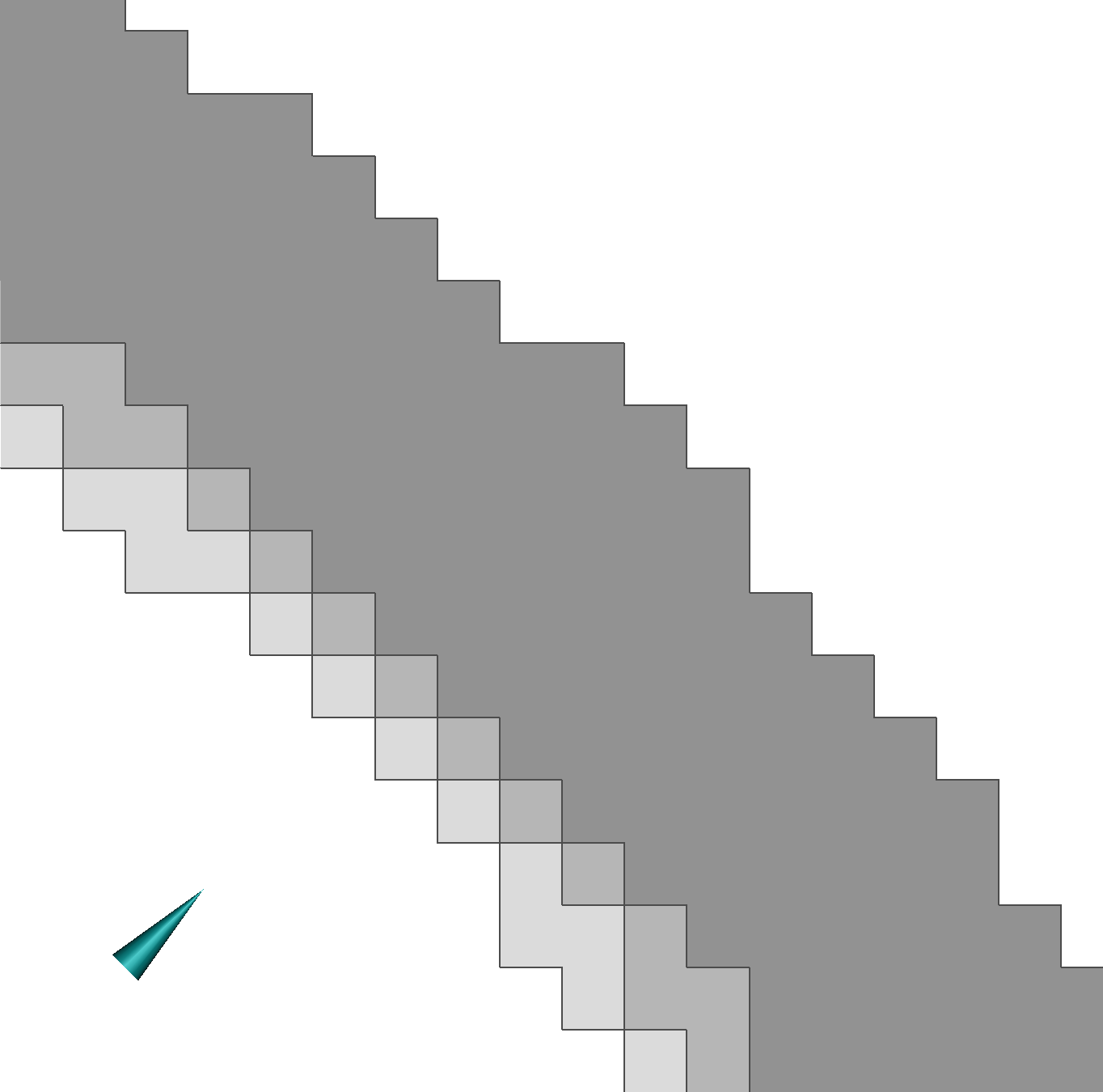} \hfill
\includegraphics[width=.24\textwidth]{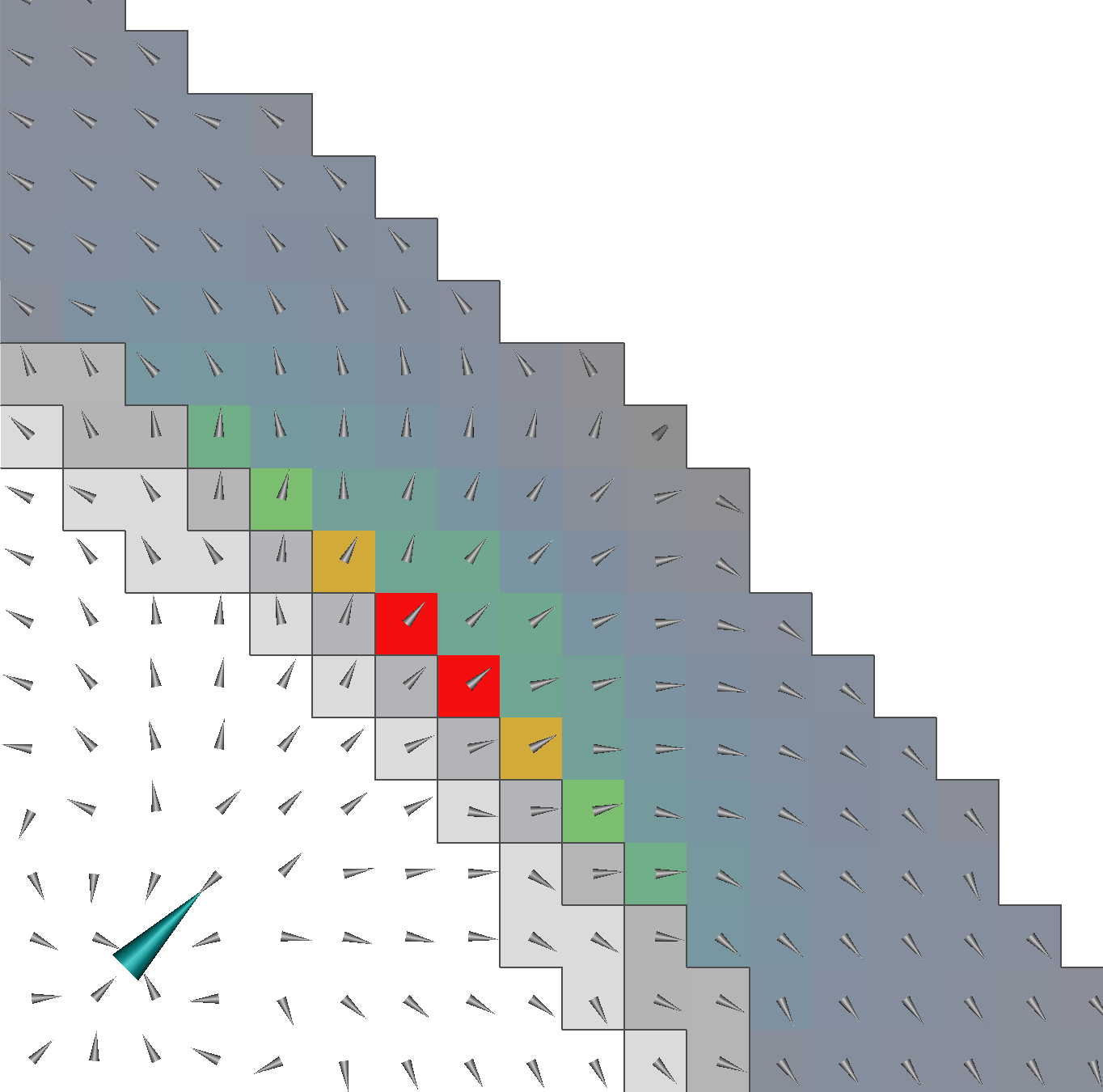} \hfil
\includegraphics[height=.22\textwidth]{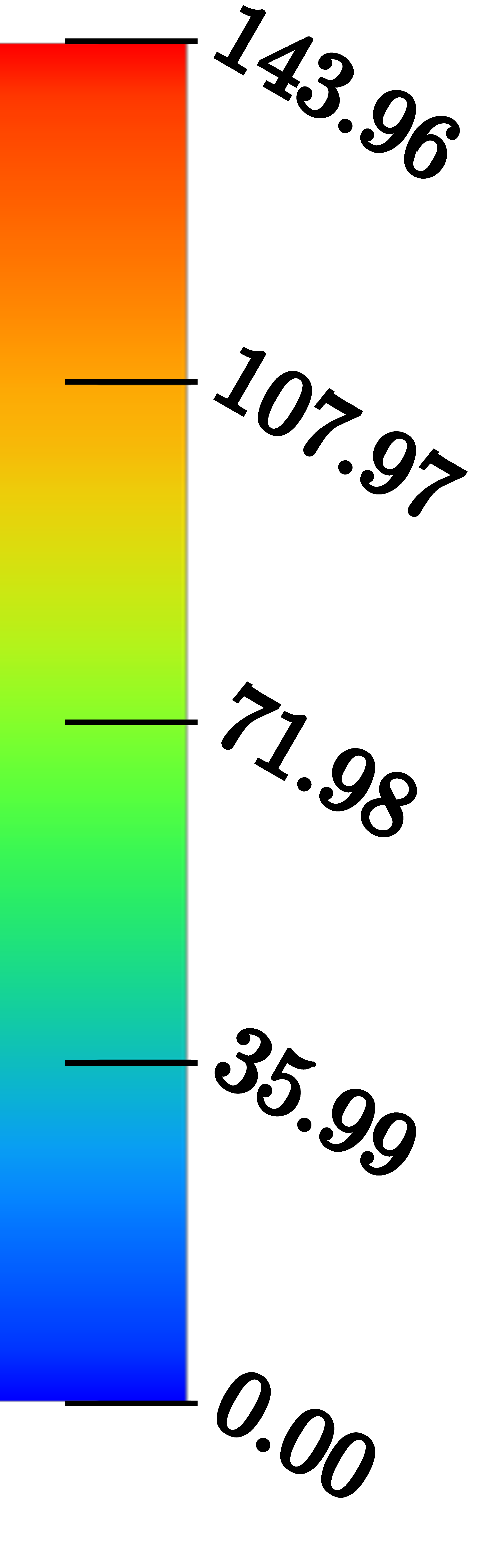} \hfill
\includegraphics[width=.24\textwidth]{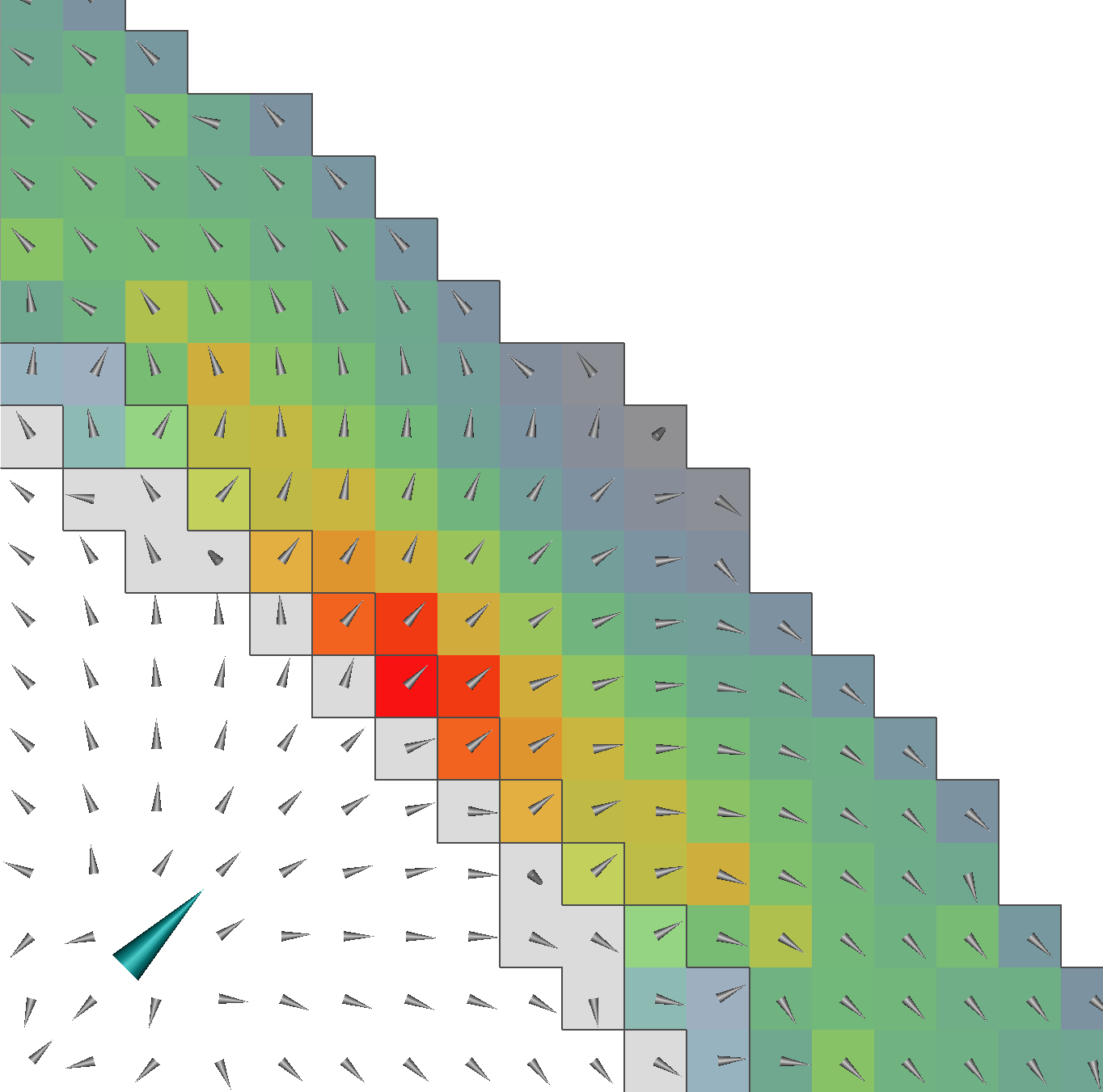} \hfil
\includegraphics[height=.22\textwidth]{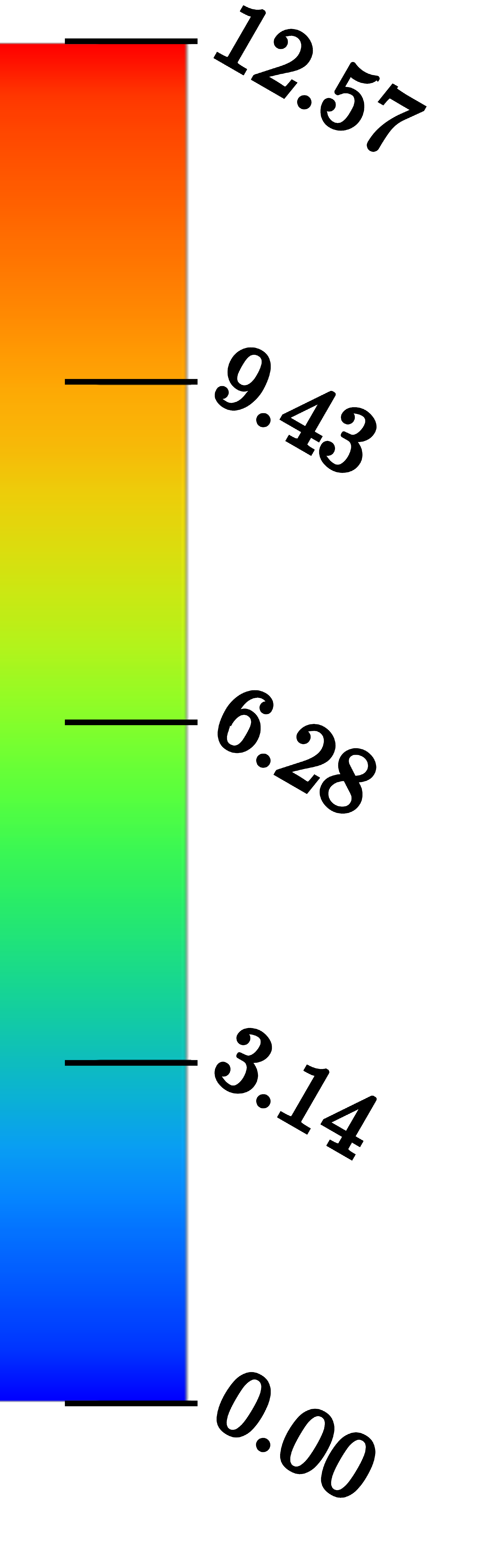}
\caption{Geometry of leaky four-layer sphere model (left, compartments from in- to outside/bottom left to top right are brain, CSF, skull, skin, and air) and visualization of strength (only skull and skin, in $\mu$A/mm$^2$) and direction of volume currents for CG-FEM (middle) and Mixed-FEM simulation (right).}%
\label{fig:current}%
\end{figure*}
The results of Sections \ref{sub:mfem-evaluation1} and \ref{sub:mfem-evaluation2} suggest that the projected Mixed-FEM is superior to the direct Mixed-FEM. To keep the presentation concise, we from here on compare only the projected Mixed-FEM with the Whitney CG-FEM.  The results for model \segresR{2}{2}{84} (Table \ref{tab:leaks}), which does not contain any skull leakages, mainly resemble those for model \segres{2}{2} for both RDM and lnMAG (Figure \ref{fig:radial-leak-mixed}).

In models \segresR{2}{2}{82} and \segresR{2}{2}{83}, the effects of the leakages become apparent. With regard to the RDM (Figure \ref{fig:radial-leak-mixed}, top row), the projected Mixed-FEM leads to lower errors in both models. In model \segresR{2}{2}{83}, the differences between the two approaches are still moderate. However, especially up to an eccentricity of 0.964 (dist. $\geq 2.8$ mm), a higher accuracy for the projected Mixed-FEM is clearly observable.
%short: For higher accuracies, the maximal errors for the projected Mixed-FEM are higher than for the Whitney CG-FEM, while median, minimal errors and IQR are lower. 
The increased number of leakages in \segresR{2}{2}{82} intensifies the difference between the approaches. The errors for the Whitney CG-FEM are clearly higher than for the Mixed-FEM here, with maximal errors larger than 0.5 at eccentricities above 0.964 (dist. $\leq 1.6$ mm).

Also with regard to the lnMAG (Figure \ref{fig:radial-leak-mixed}, bottom row), the influence of the skull leakages is apparent. In models \segresR{2}{2}{82} and \segresR{2}{2}{83}, the lnMAG increases up to an eccentricity of 0.964, and only decreases for higher eccentricities. This effect is clearly stronger for the Whitney CG-FEM than for the Mixed-FEM. In contrast, the lnMAG for the Whitney CG-FEM decreases clearly stronger than for the Mixed-FEM in model \segresR{2}{2}{84} with increasing eccentricity, leading to a switch from about 0.2 for eccentricities below 0.964 to values lower than 0.2 at an eccentricity of 0.993. Especially in model \segresR{2}{2}{83} the Whitney CG-FEM also leads to a higher variance of the lnMAG, but this variance is less distinct in the other models.

For a single, exemplary dipole, the distribution of the volume currents in skull and skin in model \segresR{2}{2}{82} simulated with the Whitney CG- and projected Mixed-FEM is visualized in Figure \ref{fig:current}. The leakage effect for the CG-FEM (Figure \ref{fig:current}, middle) is obvious. While the Mixed-FEM (Figure \ref{fig:current}, right) leads to a smooth current distribution and the highest current strengths among skull and skin elements are found in the skull compartment (up to $\approx 13$ $\mu$A/mm$^2$), the current strength peaks in the skin compartment for the Whitney CG-FEM (maximum $\approx 144$ $\mu$A/mm$^2$) and is increased by a factor of more than 11 compared to the Mixed-FEM (note the different scaling of the colorbars). Compared to the maximal current strength in the skin compartment, the current strength in the skull is very low here, showing the leakage of the volume currents through the nodes shared between CSF and the skin.

\begin{figure}[htb]
 \begin{center}
  {\includegraphics[width=0.475\textwidth]{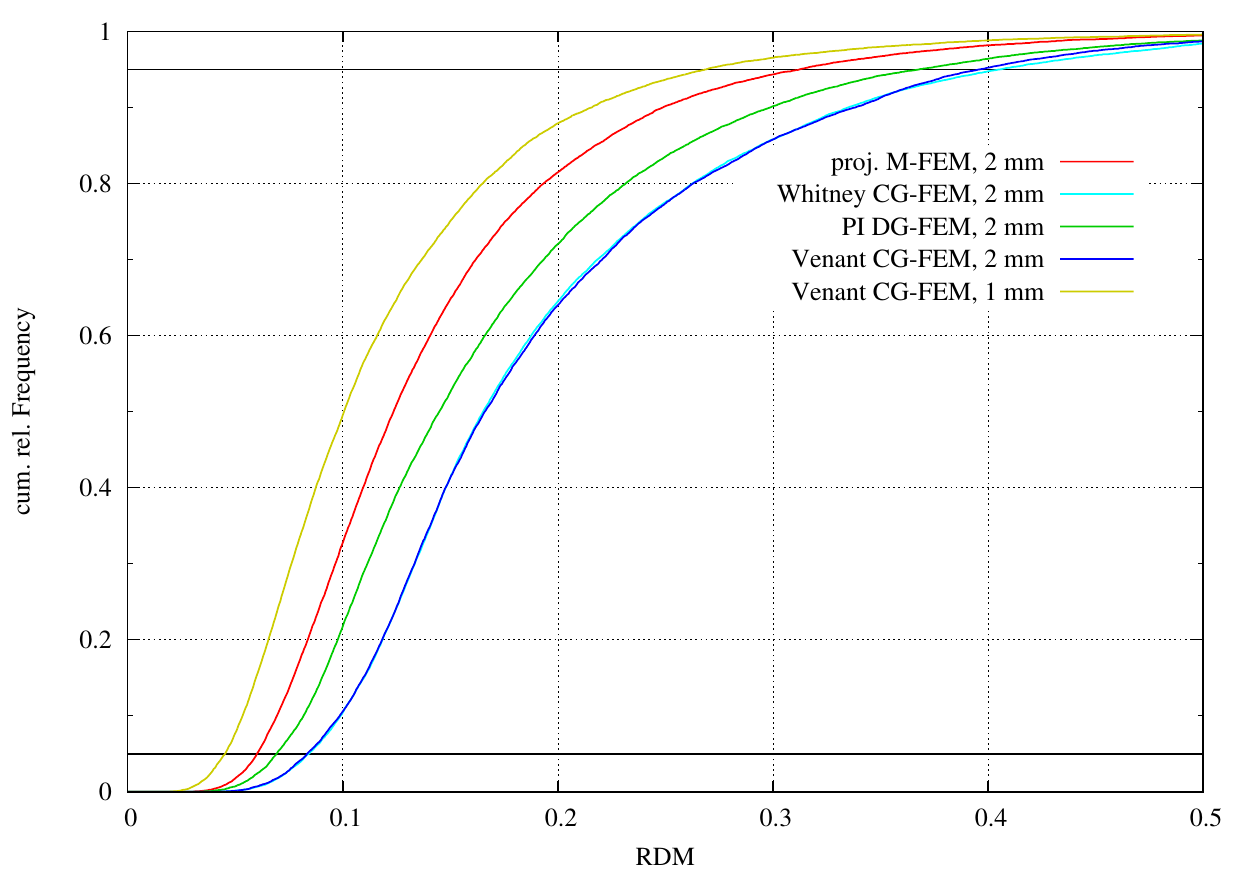}} \hfil {\includegraphics[width=0.475\textwidth]{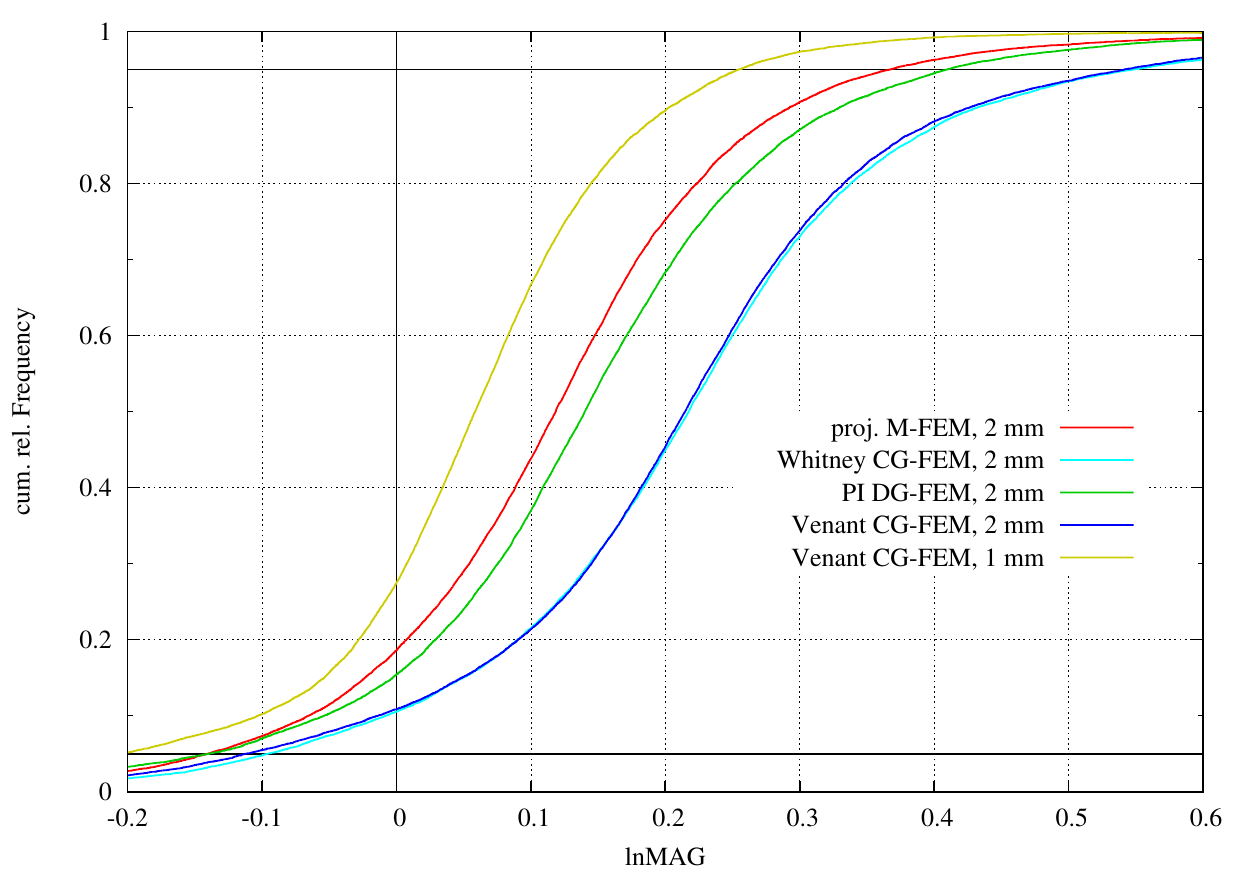}}
  \end{center}
\caption[RDM and lnMAG in realistic six-layer head model]{Cumulative relative errors of RDM (top) and lnMAG (bottom) for EEG in realistic six-layer head model. The horizontal lines indicate the 5th and 95th percentile (lower and upper lines, respectively)}
\label{fig:mixed-pi-6-realistic}
\end{figure}

\subsection{Realistic Head Model Study}
\label{sub:realistic}
The cumulative relative frequencies of RDM and lnMAG are displayed in Figure \ref{fig:mixed-pi-6-realistic}. 
Due to the rough approximation of the smooth surfaces, all models consisting of regular hexahedra (especially at the mesh width of 2 mm) lead to relatively high topography and magnitude errors when compared to the surface-based  tetrahedral reference model.
%Due to the relatively bad approximation of the geometry that is achieved when using regular hexahedra (especially at the mesh width of 2 mm), both RDM and lnMAG are relatively high.
 Comparing the results in model \ci{6}{hex}{2mm} with regard to the RDM (Figure \ref{fig:mixed-pi-6-realistic}, top), the projected Mixed-FEM performs best with roughly 95\% of the errors below 0.31
 (95\% indicated by upper horizontal bar in Figure \ref{fig:mixed-pi-6-realistic}, top).
  Therefore, the result is nearly as good as that achieved with the St. Venant approach in the 1 mm model, \ci{6}{hex}{1mm}, where 95\% of the errors are below 0.28. The partial integration DG-FEM performs nearly equally well to the Mixed-FEM with 95\% of the errors reached at about 0.36. Whitney and St. Venant CG-FEM perform nearly identically and for these approaches the 95th percentile is reached at an RDM of nearly 0.4.

With regard to the lnMAG, the differences between the results obtained using the mesh resolutions of 1 and 2 mm and also between Mixed-, DG- and the two CG-FEM approaches are larger than for the RDM (Figure \ref{fig:mixed-pi-6-realistic}, bottom). The projected Mixed-FEM performs best for model \ci{6}{hex}{2mm}, with 90\% of the errors in the range from -0.15 and 0.35 (interval between lower and upper horizontal lines in Figure \ref{fig:mixed-pi-6-realistic}). The partial integration DG-FEM performs only slightly worse with 90\% of the errors in the range from -0.15 and 0.4. Again, Whitney and St. Venant CG-FEM lead to nearly identical accuracies and show the highest errors for the model \ci{6}{hex}{2mm}, both with regard to absolute values and spread (90\% of the errors in the range from -0.1 to 0.54). The increase in accuracy when using model \ci{6}{hex}{1mm} instead of model \ci{6}{hex}{2mm} is clearer for the lnMAG than for the RDM. For the St. Venant CG-FEM, 90\% of the lnMAG-errors are in the range from -0.2 to 0.25, thus showing both a smaller spread than the results in the model \ci{6}{hex}{2mm} and also lower absolute values.

\pagebreak

\section{Discussion and Conclusion}
%\IEEEPARstart{I}{n} this paper, theoretical derivation and numerical experiments in sphere and realistic head models for an approach to solve the EEG forward problem based on the mixed finite element method (Mixed-FEM) were presented, for which existence and uniqueness of a solution as well as convergence of this solution can be proven.
\IEEEPARstart{T}{his} study introduced the Mixed-FEM approach for the EEG forward problem. 
% and a selection of a-priori error measures was stated. 
Two approaches to model the dipole source were derived, the direct and the projected. Numerical results for sphere and realistic head models were presented and compared to different established numerical methods.

The %numerical experiments show a good accuracy of the newly introduced method in common
results suggest that the Mixed-FEM achieves an appropriate accuracy for common
 sphere models, especially the projected approach. The comparison with the Whitney CG-FEM approach with optimized positions and orientations shows that the Mixed-FEM leads to comparable accuracies (Figure \ref{fig:radial-conv-mixed}). For both optimized and arbitrary source positions, the projected approach achieved a superior accuracy compared to the direct approach. 
 Previous publications concentrated on evaluating the Whitney CG-FEM in tetrahedral models \cite{bau2015comparison}. In these studies, the accuracy of the Whitney approach deteriorated when using arbitrary source positions and orientations, potentially due to the interpolation necessary to represent arbitrary source positions and orientations with the Whitney approach. This effect is not found in the hexahedral models used here and a high accuracy is achieved (Figure \ref{fig:radial-conv-mixed}). These results should be investigated in more depth in further studies. 
 In the leaky models \segresR{2}{2}{82} and \segresR{2}{2}{83}, the Mixed-FEM performs better than the Whitney CG-FEM (Figure \ref{fig:radial-leak-mixed}). This higher accuracy was expected from the Mixed-FEM based on theoretical considerations, since the Mixed-FEM is by construction charge preserving, which should prevent current leakages \cite{ewi1999superconvergence}. 
%, and also outperforms the partial integration DG approach that was introduced in the previous chapter (Figures \ref{fig:radial-geom-mixed}, \ref{fig:radial-geom-pi-mixed}).
% This is confirmed by the comparisons and visualizations in the head models with a thin skull compartment, \segresR{2}{2}{82}, \segresR{2}{2}{83}, and \segresR{2}{2}{84}. These show that, as the DG-FEM, the Mixed-FEM is leakage preventing and leads to higher accuracies and smoother -- and thereby more realistic appearing -- current distributions than the DG-FEM (Figures \ref{fig:radial-leak-pi-mixed} - \ref{fig:vis-leak-mixed}). However, it has to be kept in mind that the DG approach might lead to better results when optimizing the value of $\eta$, as discussed in Chapter \ref{ch:dg-fem}.

For EEG forward modeling, the Mixed-FEM approaches share this current preserving property with the recently proposed approaches based on the DG-FEM \cite{eng2015subtraction}. Both the direct Mixed-FEM and the partial integration DG-FEM were evaluated against CG-FEM approaches in the realistic six-compartment head model \ci{6}{hex}{2mm}. In this head model, both Mixed- and DG-FEM were advantageous in comparison to the CG-FEM (Figure \ref{fig:mixed-pi-6-realistic}). The projected Mixed-FEM clearly outperforms both Whitney and St. Venant CG-FEM in this scenario and achieves a slightly higher accuracy than the partial integration DG-FEM. Since only a few skull leakages occurred in this model and as these were concentrated in the area of the temporal bone, leakage effects do not suffice to explain the higher accuracy of Mixed- and DG-FEM. An overall higher accuracy of these approaches in this kind of model, i.e., regular hexahedral with a mesh resolution of 2 mm, can be assumed. The relatively high level of errors is a consequence of the coarse regular hexahedral meshes that were used, whereas the reference solution was computed in a highly resolved tetrahedral model. The result for the St. Venant CG-FEM in the model with a mesh resolution of 1 mm, \ci{6}{hex}{1mm}, helps to estimate the relation between the influence of the different numerical approaches and the accuracy of the approximation of the geometry. It is shown that the difference between projected Mixed-FEM and Whitney and St. Venant CG-FEM in model \ci{6}{hex}{2mm} is nearly as big as the difference between using models \ci{6}{hex}{1mm} and \ci{6}{hex}{2mm} for the St. Venant CG-FEM.

Realizing these differences in accuracy directly leads to the three main sources of error in these evaluations. Besides the previously discussed leakage effects, these are inaccurate representation of the geometry and numerical inaccuracies.
A major source of error is the representation of the geometry. Since regular hexahedral meshes were used, the influence of geometry errors is significant, especially for coarse meshes with resolutions of 2 mm or higher. No explicit convergence study comparing the results in models with increasing mesh resolution but a constant representation of the geometry was performed. However, it can be assumed from the results of previous studies that the geometry error dominates the numerical errors due to lower mesh resolutions \cite{eng2015subtraction,vor2016}.

In order to reduce the geometry error, the use of geometry-adapted meshes was considered for the CG-FEM. Such meshes have been shown to clearly improve the representation of the geometry in previous studies\mbox{\cite{CHW:Cam97,CHW:Wol2007a,wagner2016using}}. Although the use of nondegenerated parallelepipeds is uncritical for the Mixed-FEM, ``some complications may arise for general elements'' \cite{bre2012mixed}. However, it was shown that the $\Hdiv$-convergence is preserved on shape-regular asymptotically parallelepiped hexahedral meshes\mbox{\cite{ber2005approximation}} and, for the two-dimensional case, error estimates for general quadrilateral grids can be obtained when modifying the lowest-order Raviart-Thomas elements \cite{cho2002flux,kwa2011mixed} and for convex quadrilaterals even superconvergence was shown \cite{ewi1999superconvergence}. The use of geometry-adapted hexahedral meshes in combination with the Mixed-FEM should therefore be evaluated in future studies.

Regarding the numerical inaccuracy due to the discretization of the equations and the source singularity, the Mixed-FEM allows to increase the regularity of the right-hand side by one degree. As a consequence of the first\replaced{-}{ }order formulation \eqref{eq:mixed-weak}, applying the derivative to the delta distribution included in the primary current $\j^p$ can be circumvented. The results obtained show high numerical accuracies, especially at the highest eccentricities, and particularly for the projected Mixed-FEM. 
This increase in accuracy comes at the cost of a higher number of degrees of freedom than that of the CG-FEM, as the current $\j$ is also considered as an unknown now, meaning that it has to be discretized. Furthermore, the discrete problem has a saddle point structure \eqref{eq:mixed-disc} and cannot be efficiently solved with AMG-CG solvers without further modifications. 
Although the number of unknowns is clearly increased compared to the CG-FEM, e.g., in model \segres{2}{2} we have $\text{\#DOF}_{M} = 1,243,716 + 407,904$, and $\text{\#DOF}_{CG} = 428,185$ (cf. Table \ref{tab:models}), by introducing an algorithm based on the idea of the conjugated Uzawa-iteration (Section \ref{sub:solving}), the solving time even in the finest model \segres{1}{1} was reduced to less than two minutes. This solving time is only a few seconds slower than that for the CG-FEM. Furthermore, as the equation system \eqref{eq:mixed-disc} is symmetric, the transfer matrix approach \cite{CHW:Wei2000, CHW:Wol2004} can be applied for the Mixed-FEM to reduce the number of equation systems that have to be solved to equal the number of sensors.

As an alternative to the straightforward approach presented here for solving the linear equation system \eqref{eq:mixed-disc-system} using the Schur complement, an approach based on the method of Lagrange multipliers has been proposed \cite{de1968equilibrium}. In this approach, the continuity of the vector-valued basis functions is no longer enforced by the definition of the basis functions, but by introducing inter\deleted{-}element Lagrange multipliers. This approach leads to a linear equation system having as many unknowns as the number of faces in the case of lowest-order Raviart-Thomas elements. This equation system is symmetric, positive definite, and sparse. Although this approach does not necessarily lead to a decrease of the solving time \cite{ver1984combined,ber1994mixed}, a higher order of convergence is predicted in theory when employing the information contained in the Lagrangian multipliers \cite{arn1985mixed,bre2012mixed}.
 Therefore, it is desirable to evaluate this solution approach in subsequent studies.
 
The lowest-order Raviart-Thomas elements used in this study are the most classical, but only one of many different elements that have been developed to approximate $\Hdiv$. Further element types are, e.g., Brezzi-Douglas-Marini (BDM)\mbox{\cite{bre1985two,bre1987mixed}} and Brezzi-Douglas-Fortin-Marini (BDFM)\mbox{\cite{bre1987efficient}} elements. To overcome known limitations of these classical element types, further elements to approximate $\Hdiv$ were developed more recently\mbox{\cite{arn2005quadrilateral,fal2011hexahedral}}. Due to different approximation properties of the element types, the evaluation of further element types for solving the EEG forward problem using the Mixed-FEM in future studies might be worthwhile. Also the use of higher-order Raviart-Thomas elements, e.g., $RT_1$ elements in combination with discontinuous linear Ansatz-functions for the potential, should be considered, as the theoretically predicted convergence rates improve for higher element orders. For an overview of the most common finite element spaces to approximate $\Hdiv$, including higher-order elements, and their convergence properties, we refer the reader to\mbox{\cite{bre2012mixed}}. However, the use of higher order elements comes at the cost of an increased number of degrees of freedom. Thus, the use of higher mesh resolutions should always be considered as an alternative to the use of higher-order elements.
%
%Especially in the context of the high errors due to the inexact representation of the geometry, also the use of elements that are specifically tailored for non-regular hexahedrons, as they occur when applying geometry-adaptation techniques, should be considered.

%Therefore, future experiments should investigate the accuracy of the cut-cell DG approach in comparison to the Mixed-FEM using regular or geometry-adapted hexahedral meshes for which the convexity of the elements is ensured. Furthermore, the accuracies in realistic head models with finer mesh widths than the here used 2 mm should be evaluated. This is uncritical for Mixed- and CG-FEM, while the memory demand for the DG-FEM is in this case in a range that would be too high for nowadays common PCs.

As mentioned, the Mixed-FEM guarantees the conservation of charge by construction. In consequence, especially in models with thin insulating compartments and at highest eccentricities, it still leads to high accuracies, which also encourages the use of the Mixed-FEM in related applications that depend on an accurate simulation of the electric current, such as the magnetoencephalography (MEG) forward problem, transcranial direct current stimulation (tDCS), or deep brain stimulation (DBS) simulations.

%short? nach abstract checken!
% For MEG, the simulation of the secondary magnetic flux $\Phi_s$ directly depends on an accurate simulation of the volume currents and it was shown that this can be locally inaccurate when using CG-FEM approaches, especially when leakages occur in the model (Figure \ref{fig:current}). However, the effects of leakages are rather local and the currents incorrectly penetrating the skull are relatively weak compared to those in the vicinity of the source position. This aside, the MEG solution should also in non-leakage scenarios profit from a physically correct and charge preserving simulation of the volume currents, as it is guaranteed by the Mixed-FEM. Thus, the use of Mixed-FEM for MEG forward simulations should be considered.

Overall, we conclude that the Mixed-FEM is an interesting new approach that can at least complement and in some scenarios even outperform standard continuous Galerkin FEM approaches for simulation studies in bioelectromagnetism. The use of different element types and solving algorithms should be investigated in further studies. 

% or


\nopagebreak

\section*{Acknowledgment}
The authors would like to thank the anonymous reviewers for their valuable comments and suggestions to improve the quality of the paper. We are also grateful to Prof. Dr. Steffen B\"orm for proof-reading and his advice with regard to the Theory section.

\ifCLASSOPTIONcaptionsoff
  \newpage
\fi



%\bibliographystyle{IEEEtran}

\bibliographystyle{IEEEtran}
\bibliography{IEEEabrv,./mixed}

\end{document}